\documentclass[12pt]{article}

\usepackage[latin1]{inputenc}
\usepackage{amsmath,amsthm,amssymb,amscd,}
 \newtheorem{thm}{Theorem}[section]
 \newtheorem{cor}[thm]{Corollary}
 \newtheorem{lem}[thm]{Lemma}
 \newtheorem{prop}[thm]{Proposition}
 \theoremstyle{definition}
 \newtheorem{df}[thm]{Definition}
 \theoremstyle{remark}
 \newtheorem{rem}[thm]{Remark}
 
 \numberwithin{equation}{section}

\def\be#1 {\begin{equation} \label{#1}}
\newcommand{\ee}{\end{equation}}
\def\dem {\noindent {\bf Proof : }}

\def\sqw{\hbox{\rlap{\leavevmode\raise.3ex\hbox{$\sqcap$}}$%
\sqcup$}}
\def\findem{\ifmmode\sqw\else{\ifhmode\unskip\fi\nobreak\hfil
\penalty50\hskip1em\null\nobreak\hfil\sqw
\parfillskip=0pt\finalhyphendemerits=0\endgraf}\fi}

\newcommand{\mb}{\medskip\noindent}
\newcommand{\gb}{\bigskip\noindent}
\newcommand{\R}{\mathbb R}

\newcommand{\N}{\mathbb N}
\newcommand{\Z}{\mathbb Z}

\newcommand{\s}{\mathcal S}

\title{Uniform estimates for paraproducts and related multilinear multipliers.}
\author{Fr\'ed\'eric Bernicot \\
\small Universit\'e de Paris-Sud, Orsay et CNRS 8628, \\ \small 91405 Orsay Cedex, France \\
\small {\em E-mail address:} {Frederic.Bernicot@math.u-psud.fr} }

\begin{document}
\maketitle

\begin{abstract}
 In this paper, we prove some uniform estimates between Lebesgue and Hardy spaces for operators closely related to the multilinear paraproducts on $\R^d$. We are looking for uniformity with respect to parameters, which allow us to disturb the geometry and the metric on $\R^d$.
\end{abstract}

\section{Introduction}

The purpose of this article is to prove uniform estimates on paraproducts and similar multilinear operators. Let us first recall what is a paraproduct. \\
A $n$-linear paraproduct $\Pi$ on $\R^d$ is a $n$-linear operator of the following form~:
$$ \Pi(f_1,..,f_n)(x):= \int_0^\infty \prod_{i=1}^{n} \pi^i_t \ast f_i(x) \frac{dt}{t}, $$
or of the discrete form
$$ \Pi(f_1,..,f_n)(x):= \sum_{j\in\Z} \prod_{i=1}^{n} \pi^i_{2^j} \ast f_i(x) \frac{dt}{t}.$$
Here the $\pi^i_t$ are smooth functions which Fourier transform $\widehat{\pi^i_t}$ are bump functions adapted to the ball $\{ \xi\in \R^d,\ |\xi|\leq 1\}$ and we assume that there exists one index $i\in\{1,..,n\}$ such that
$$\forall t>0, \qquad   \widehat{\pi^i_t}(0).$$
In all the sequel, a smooth function $\phi$ is said to be ``adapted to a set'' $I\subset \R^d$ if it is supported on this set and satisfies : for all order $n\in \N^d$
$$\|\phi^{(n)} \|_\infty \leq |I|^{-|\alpha|}.$$
Then for such a paraproduct, the classical Calder\'on-Zygmund theory gives us that for all exponents $1<p_1,..,p_n<\infty$ such that
$$ 0<\frac{1}{p}:= \sum_{i=1}^{n} \frac{1}{p_i}<1,$$
there exists a constant $C=C(p_i)$ such that for all functions $f_i\in \s(\R^d)$,
$$ \left\| \Pi(f_1,..,f_n) \right\|_{p} \leq C \prod_{i=1}^{n} \|f_i\|_{p_i}.$$ 
These estimates in Lebesgue spaces depend on the functions $\pi^i_t$. We would like to understand how can we modify these functions, keeping uniform estimates. \\

\mb The paraproducts are the first studied class of singular bilinear operators. Their study began by the works of  J.M. Bony in \cite{bony} and of R. Coifman and Y. Meyer in \cite{CM2,cm,CM3}, where in particular continuities in Lebesgue spaces are shown. The first uniform result is the following one (from \cite{CM4})~:
\begin{thm} \label{thm1} Let $(M_i)_{1\leq i\leq n}$ be integers and $(\pi_{j}^{i})_{\genfrac{}{}{0pt}{}{j\in\Z}{1\leq i\leq n}}$ be smooth functions such that $\widehat{\pi_{j}^{i}}$ is adapted to the rectangle $[-2^{j+M_i},2^{j+M_i}]^{d}$. Assume that there exist an integer $N$ and an index $i\in\{1,..,n\}$ such that for all $j$,
\be{hyp1}  \forall \eta \in [-2^{j+M_i-N},2^{j+M_i-N}]^{d}, \qquad \widehat{\pi_{j}^{i}}(\eta)=0 . \ee
Then for all exponents $1<p_1,...,p_n\leq\infty$ satisfying
$$0< \frac{1}{p}:=\frac{1}{p_1}+..+\frac{1}{p_n}<1,$$
there exists a constant $C=C(N,p_i)$, which does not depend on $(M_i)_i$ such that
$$\forall f_i\in\s(\R^d), \qquad \left\| \sum_{j} \prod_{i=1}^{n} \pi_{j}^{i} \ast f_i \right\|_{p} \leq C \prod_{i=1}^{n} \|f_i\|_{p_i}.$$
\end{thm}

\mb About this result, there are two different questions : what is the maximal range of exponents with uniform estimates ? May we weaken the assumption (\ref{hyp1}) ?

\mb The second question was solved by C. Muscalu, T. Tao and C. Thiele in \cite{MTT4}, where they prove the stronger result~:

\begin{thm} \label{thm2} Let $(M_i)_{1\leq i\leq n}$ be integers and $\pi_{j}^{i}$ be smooth functions such that $\widehat{\pi_{j}^{i}}$ be adapted to the rectangle $[-2^{j+M_i},2^{j+M_i}]^{d}$. Assume that for all $j$, there exists an index $i\in\{1,..,n\}$ with
\be{hyp2} \widehat{\pi_{j}^{i}}(0)=0. \ee
Then for all exponents $1<p_1,...,p_n<\infty$ satisfying
$$0<\frac{1}{p}:=\frac{1}{p_1}+..+\frac{1}{p_n}<1,$$
there exists a constant $C=C(N,p_i)$, which does not depend on $(M_i)_i$ such that
$$\forall f_i\in\s(\R^d), \qquad \left\| \sum_{j} \prod_{i=1}^{n} \pi_{j}^{i} \ast f_i \right\|_{p} \leq C \prod_{i=1}^{n} \|f_i\|_{p_i}.$$
\end{thm}

\mb It is even shown a little stronger version (a maximal version) than this one. The assumption (\ref{hyp2}) is much weaker than (\ref{hyp1}). The proof of Theorem \ref{thm2} is a mixture of the proof of Theorem \ref{thm1} and arguments from graph theory. \\
Such a result was motivated by the paper \cite{mtt} from the same authors. In this article, they study some uniform estimates for multilinear operators far more singular than the paraproducts, closely related to the bilinear Hilbert transforms. The ``classical'' time-frequency analysis, to decompose these kind of operators, uses some information and estimates on operators, which look like paraproducts. That is why they have first shown in \cite{MTT4} uniform estimates for paraproducts. 

\mb In this paper, we are interested in answering to the first question. Mainly we want to obtain uniform estimates with infinite exponents and some exponents lower than $1$. In \cite{li}, X. Li has shown uniform estimates when $1<p_1,...,p_n<\infty$ and $p$ may be lower than one. We would like to extend his result for some exponent $p_i<1$ or $p_i=\infty$.
The continuities for this range of exponents have already been proved (for example in \cite{GK} by L. Grafakos and N. Kalton). Here we would like to improve these continuities with uniform estimates. \\
A second motivation for the study of paraproducts is this one : we know how to decompose a multilinear multiplier, satisfying H\"ormander's condition, with multilinear paraproducts. A $n$-linear multiplier $T$ is given by its symbol $\sigma \in \s(\R^{dn})$, with the formula~:
\be{opsym} T(f_1,..,f_n)(x):= \int_{\R^{dn}} e^{ix.(\xi_1+..+\xi_n)} \sigma(\xi) \prod_{i=1}^n \widehat{f_i}(\xi_i) d\xi. \ee
The H\"ormander condition corresponds to the following assumption~:
\be{marcin} \forall m_i\in\N^d, \qquad \left| \prod_{i=1}^{n} \partial_{\xi_i}^{m_i} \sigma(\xi_1,..,\xi_n) \right| \leq \frac{1}{\left(|\xi_1|+..+|\xi_n|\right)^{|m_1|+..|m_n|}}. \ee
Note the appearance of the quantity $|\xi_1|+..+|\xi_n|$, which corresponds to the distance $d(\xi,0)$ in the frequency plane. We are now interested in disturbing the metric. We would like study the following distance
$$ d_\lambda(\xi,0):=\sum_{i=1}^{n} |\lambda_i \xi_i|,$$
given by non vanishing reals $\lambda_i$. In fact it is easy to see that our parameters $\lambda_i$ have the same function than the parameters $M_i$ of Theorems \ref{thm1} and \ref{thm2} (we have the relation $\lambda_i \simeq 2^{-M_i}$). So we would like to have uniform estimates with respect to the new distance $d_\lambda$.
The problem of disturbing the metric appeared for example in the study of bilinear Hilbert transforms along polynomial curves (\cite{lifan}) and was one of the X. Li motivations to study uniform estimates for paraproducts.

\mb We will also prove the following result~:

\begin{thm} \label{thm:general} Let $\sigma$ be an $x$-independent symbol such that
\be{marcin1} \forall m_i \in \N^{d} \qquad \left|\partial_{\xi_1}^{m_1} .. \partial_{\xi_{n-1}}^{m_{n}} \sigma(\xi_1,..,\xi_{n}) \right| \lesssim \frac{\prod_{i=1}^{n} |\lambda_i|^{|m_i|}} {d_\lambda(\xi,0)^{|m_1|+..+|m_{n-1}|}}. \ee
Let $0<p_i,p\leq \infty$ exponents satisfying
$$\frac{1}{p}=\sum_{i=1}^{n} \frac{1}{p_i}.$$
Let us denote the three disjoint sets (which may be empty) $S_1,S_2$ and $S_3$ such that
$$\{1,..,n\}=S_1 \sqcup S_2 \sqcup S_3,$$
with
$$ \forall i\in S_1,\ p_i=1, \qquad \forall i\in S_2,\ p_i=\infty \quad \textrm{and} \quad \forall i\in S_3,\ p_i\in\{1,\infty\}^c.$$
Then we know that the multilinear multiplier  $T$ defined by (\ref{opsym}) can be continuously extended from $\otimes_{i=1}^{n} F_i$ to $G$ in the three following cases~: 
\begin{itemize}
 \item if $0<p<\infty$ with $G=L^p$, $F_i=H^{p_i}$ for $i\in S_1\cup S_3$ and $F_i=L^\infty_c$ for $i\in S_2$,  
 \item if $0<p<\infty$ with $G=L^{p,\infty}$, $F_i=L^1$ for $i\in S_1$, $F_i=H^{p_i}$ for $i\in S_3$  and $F_i=L^\infty_c$ for $i\in S_2$,
 \item if $p=\infty$ (and also for all $i\in\{1,..,n\}$ $p_i=\infty$) with $G=BMO$ and $F_i=L^\infty_c$ for all $i\in\{1,..n\}$. 
\end{itemize}
In addition we have the two following improvements~: \\
Part 1) : All these continuities are uniformly bounded with respect to the parameters $\lambda_i$ under one of the two following assumptions~:
\begin{align*}
 a-) & \qquad \forall 1\leq i\leq n, \qquad p_i<\infty \\
 b-) & \qquad \displaystyle \gamma:=\sum_{\genfrac{}{}{0pt}{}{j \in\{1,..,n\} }{ |\lambda_j| \simeq \max\{ |\lambda_l|,\ 1\leq l \leq n\} }} \frac{1}{p_j} \geq \frac{1}{2}.
  \end{align*}
Else the continuity bound depends on the ratio
$$\frac{\max\{ |\lambda_k|,\ 1\leq k \leq n\}}{\min\{ |\lambda_k|,\ 1\leq k \leq n\}}.$$
Part 2) : We don't know if the conditions a-) or b-) are sufficient to get uniform bounds, however we will show that if $p<\infty$ and $\gamma=0$ then we cannot have a uniform bound.
\end{thm}

\mb In this result we write $L^p=L^p(\R^d)$ for the ``classical'' Lebesgue spaces and $H^p=H^p(\R^d)$ for the Hardy spaces (which is equivalent to the Lebesgue spaces $L^p$ if $1<p<\infty$) and $BMO=BMO(\R^d)$ the space of functions of ``bounded mean oscilaation''. We write $L^\infty_c$ for the set of bounded compactly supported functions, equipped with the $L^\infty$-norm. 

\begin{rem} Our proof, show that in particular case, we can obtain the continuity with the whole space $L^\infty$ instead of $L^\infty_0$. For convenience and technical difficulties (see the proof of Corollary \ref{cor2}), we prefer to only work with the space $L^\infty_0$.
\end{rem}

\begin{rem} By Taking $\lambda_i \simeq 2^{-M_i}$, the paraproducts of Theorems \ref{thm1} and \ref{thm2} verify (\ref{marcin1}) uniformly with respect to $\lambda$ because the symbol is given by
 $$ \sigma(\xi_1,..,\xi_n)=\sum_j \prod_{i=1}^n \widehat{\zeta_{i,j}}(\xi_i).$$
So Theorem \ref{thm:general} improves the uniform estimates of Theorem \ref{thm2} and answer to the asked question.
\end{rem}

\begin{rem} By using time-frequency tools such as ``tiles'' and ``trees'' as in \cite{mptt,li}, it should be possible to prove some uniform ``weak type restricted estimates'' in $L^p$, which are stronger than our continuity in $L^p$ for $p<1$.
\end{rem}

\begin{rem} The continuities are already known from the papers \cite{GK} and \cite{GK2} of L. Grafakos and N. Kalton. In fact our operators are multilinear Calder\'on-Zygmund operators and so their continuities are a consequence of the paper \cite{GT} of L. Grafakos and R. Torres. The improvement is the fact that we can have uniform bounds and we must be careful because the constants, as multilinear Calder\'on-Zygmund operators, are not uniformly bounded. So we will use the ideas of the Calder\'on-Zygmund theory with a few improvements.
\end{rem}

\mb There is an other interest to study such uniform estimates. The symbols verifying (\ref{marcin1}) uniformly with respect to $\lambda$ satisfy the Marcinkiewicz condition~:
\be{marcin2} \forall m_i \in \N^{d} \qquad \left|\partial_{\xi_1}^{m_1} .. \partial_{\xi_{n-1}}^{m_{n}} \sigma(\xi_1,..,\xi_{n}) \right| \lesssim \prod_{i=1}^{n} |\xi_i|^{-|m_i|} . \ee 
However, from \cite{GK} we know that the condition (\ref{marcin2}) is in general not sufficient to guarantee continuity, as in the previous Theorem. So our result allows us to almost describe the ``limit case'' between (\ref{marcin}) and (\ref{marcin1}) to get these continuities.

\mb To prove our Theorem, we will use model operators, which generalize and are more symmetric than the paraproducts. In the definition of paraproducts, there has to be one (or more) index $i\in\{1,..,n\}$ such that (\ref{hyp1}) or (\ref{hyp2}) is satisfied, so there is a lack of symmetry in their definition (see Remark \ref{remsym}). 

\gb The plan of this paper is the following one. In Section \ref{notations}, we define notations and our model operators. We first prove Theorem \ref{thm:general} for our model operators~: in the case where all exponents belong to $(1,\infty)$ in Section \ref{section:littlewood} (this part only uses Littlewood-Paley theory) and after for others exponents in Section \ref{section:carleson} (this part uses Carleson measures and an improved Calderon-Zygmund theory). Then we complete the proof of Theorem \ref{thm:general} for general multipliers in Section \ref{section:paradecomp}.

\section{Definition of our model operators.}
\label{notations}

For the rest of this paper, we use the well-known notations~:
let $\zeta$ be a function on $\R^d$, $t\neq 0$ be a real and $q\in\R^d$ be a vector. We set $\zeta_t$ and $\zeta_{t,q}$ for the $L^1$-normalized functions defined by
$$\zeta_{t}(x) := \frac{1}{|t|^d} \zeta(t^{-1}x) \quad \textrm{and} \quad \zeta_{t,q}(x) := \frac{1}{|t|^d} \zeta(t^{-1}(x-q)).$$

\mb We will work with the Hardy spaces on $\R^d$, so let us first recall one of its definitions.

\begin{df} Let $\Psi$ be a smooth function. We define $S_\Psi$ to be the continuous or the discrete Littlewood-Paley square function, given by~
$$S_\Psi(f):=\left( \int \left| \Psi_t \ast f\right|^2 \frac{dt}{t} \right)^{1/2} \quad \textrm{or} \quad S_{\Psi}^\Delta(f):=\left( \sum_{n\in\Z} \left| \Psi_{2^n}
\ast f\right|^2 \right)^{1/2}.$$
\end{df}

\mb We use these functionals to get the following definition of Hardy spaces (See \cite{Gra})~:
\begin{df} \label{hardyspace} Let $\Psi$ be a non null smooth function whose spectrum is contained in a corona around $0$. For $0<p<\infty$, we define the Hardy space $H^p=H^p(\R^d)$ as the set of distributions $f\in \s'(\R^d)$ satisfying~:
$$ \|f\|_{H^p}:=\|S_\Psi (f)\|_{p}<\infty.$$
From the book \cite{stein} we know that for $1<p<\infty$ the Hardy space $H^p$ corresponds to the Lebesgue space  $L^p$. In addition, we have the choice to keep a discrete or a continuous square function : the definition of the space does not depend on it or on the choice of the function $\Psi$.
\end{df}

\mb We have to control norms in the Schwartz space, so we set for an integer $K$
\be{coeffck} c_K(\zeta) := \sup_{x\in \R^d} \left(1+|x| \right)^{K} \sup_{\genfrac{}{}{0pt}{}{\alpha\in \N^d}{|\alpha|\leq K}} \left| \partial_x^\alpha \zeta(x) \right|. \ee

\mb Now we define our model operators.

\begin{df} Let $\Psi$ be a smooth function on $\R^d$ whose spectrum is contained in a corona around $0$ and let $\Phi^i$ be smooth functions whose spectrum is bounded. Let $L$ be a bounded function on $\Z$, $\lambda=(\lambda_1,..,\lambda_{n})\in(\R^*)^{n}$ and $\rho=(\rho_1,..,\rho_{n})\in ]0,1]^{n}$ be parameters. Then we define the following operator~:
$$T_{\rho,\lambda,L}(f_1,..,f_{n})(x)= \sum_{k\in\Z}
 L(k) \int_{\R^d} \Psi_{2^k}(y)\prod_{i=1}^{n}
 [\Phi^{i}_{\lambda_i 2^k}\ast f_i](x-\rho_i\lambda_i y)dy.$$
We also have the continuous version with a bounded function $L$ on $\R^{+}$, defined by
$$U_{\rho,\lambda,L}(f_1,..,f_{n})(x)= \int_{0}^\infty
 L(t) \int_{\R^d} \Psi_{t}(y)\prod_{i=1}^{n}
 [\Phi^{i}_{\lambda_i t}\ast f_i](x-\rho_i\lambda_i y)dy \frac{dt}{t}.$$
\end{df}

\mb It is easy to see that these operators continuously act from $\s(\R^d)^{\otimes n}$ to $\s(\R^d)$. In addition, the operator $T_{\rho,\lambda,L}$ is associated to the following symbol $\sigma$~:
$$ \sigma(\xi_1,..,\xi_n):= \sum_{k\in\Z}
 L(k) \widehat{\Psi}\left(2^k(\rho_1\lambda_1\xi_1+..+\rho_n\lambda_n\xi_n )\right)\prod_{i=1}^{n}
 \widehat{\Phi^{i}} (\lambda_i 2^k \xi_i),
$$
which satisfies (\ref{marcin1}) uniformly on $\lambda$.

\mb We want now to make the link with the ``classical'' paraproducts.

\begin{prop} \label{proppara} The parameters $\rho_i$ allow us to get the ``classical'' paraproducts as limit of our previous operators~: for all $f_1,..,f_n\in\s(\R^d)$
\be{limit} U_{\rho,\lambda,L}(f_1,..,f_{n})(x)
\xrightarrow [] { \genfrac{}{}{0pt}{}{\rho_1=1}{\rho_i \to 0} } \int_{0}^\infty
L(t)  [\Psi \ast \Phi^{1}]_{\lambda_1 t}\ast f_1(x) \prod_{i=2}^{n}
 [\Phi^{i}_{\lambda_i t}\ast f_i](x) \frac{dt}{t}.\ee
Here the convergence is in the $\s(\R^d)$ sense.
\end{prop}

\mb We do not write the details of this result. With the good assumptions about the functions $f_i$, it is easy to prove this convergence. 

\begin{rem} \label{remsym} Our model operators have a symmetry : the definition is invariant by permutations on the $n$ functions, which is not the case for the ``classical'' paraproducts. For example in the bilinear case, we want to estimate in $L^2$ the two different paraproducts (for $f\in L^\infty$ and $g\in L^2$)~:
$$ \int_0^\infty \left[\Psi_{\lambda_1 t}\ast f \right] \left[ \Phi_{\lambda_2 t} \ast g \right] \frac{L(t)}{t} dt \qquad \textrm{and} \qquad  \int_0^\infty \left[\Phi_{\lambda_1 t}\ast f\right] \left[\Psi_{\lambda_2 t} \ast g \right] \frac{L(t)}{t} dt, $$
uniformly on $(\lambda_1,\lambda_2)$ with $|\lambda_2|>|\lambda_1|$. These two paraproducts are a little different and so their study ask some different arguments.  \\
That is why we prefer working with our model operators, which own symmetry invariance and allow us to get by a limit argument these two kinds of paraproducts.
\end{rem}

\begin{rem} It is quite easy to show that our model operators satisfy the assumptions of Theorem \ref{thm:general} with uniform bounds with respect to $\lambda$ and $\rho$. We let to the reader the details of this claim.
\end{rem}

\mb Before to prove the positive part (part 1) of Theorem \ref{thm:general} for our model operators, we would like to explain the negative claim of this Theorem (part 2) in the bilinear case~:

\begin{prop} Let $\rho_1=\rho_2=1$ be fixed and $|\lambda_1| \ll |\lambda_2|$ be reals. There exists operators $U_{\rho,\lambda,m}$ (also satisfying the assumptions of Theorem \ref{thm:general}) which cannot be continuous from $L^p\times L^{\infty}$ into $L^p$ for $1<p<\infty$ with an uniform bound with respect to $\lambda$. 
\end{prop}

\dem Let us choose $\Phi^{i}=\zeta$ a smooth and nonnegative function whose integral is equal to $1$ and set
$$ U_{\epsilon,\lambda}(f,g)(x):= \int_\epsilon^{\epsilon^{-1}}\int_{\R^d} \Psi_{t}(y) \zeta_{\lambda_1 t}\ast
f(x-\lambda_1 y)\zeta_{\lambda_2t}\ast g(x-\lambda_2y) \frac{dydt}{t}.$$ When  $\lambda_1$ tends to $0$, we have
$$\forall f\in \s(\R^d),\, x\in\R^d \qquad  \lim_{\lambda_1\rightarrow 0} \zeta_{\lambda_1 t}\ast
f(x-\lambda_1 y) = f(x).$$
Due to the presence of the $\epsilon>0$, we have for $f,g\in\s(\R^d)$ 
$$\forall x\in\R^d \qquad \lim_{\lambda_1\rightarrow 0} U_{\epsilon, \lambda}(f,g)(x) = f(x)\int_{\epsilon}^{\epsilon^{-1}} \int_{\R^d} \Psi_t(y) \zeta_{\lambda_2 t}\ast g(x-\lambda_2 y)
\frac{dy dt}{t}.$$
We can now take $\epsilon\to 0$ and we get
$$ \lim_{\epsilon \to 0} \int_{\epsilon}^{\epsilon^{-1}} \int_{\R^d} \Psi_t(y) \zeta_{\lambda_2 t}\ast g(x-\lambda_2 y)\frac{dy dt}{t} = \int_{0}^{\infty} \int_{\R^d} \Psi_t(y) \zeta_{\lambda_2 t}\ast g(x-\lambda_2 y)
\frac{dy dt}{t}.$$
Then we can choose good functions $\Psi$ and $\zeta$ in order to find the linear Hilbert transform $H$. With these ones, we conclude 
$$\forall x\in\R^d, \qquad \lim_{\epsilon \to 0} \lim_{\lambda_1\rightarrow 0} U_{\epsilon,\lambda}(f,g)(x) = f(x)H(g)(x).$$
So if we have uniform estimates on $U_{\epsilon, \lambda}$ from $L^p \times L^\infty$ into $L^p$, by using Fatou's lemma, we get~:
$$ \forall f,g\in \s(\R^d), \qquad \left\| fH(g)\right\|_{p} \lesssim \|f\|_{p}\|g\|_\infty.$$
Such an estimate implies the boundedness of $H$ on $L^\infty$ which is not possible. So we cannot have uniform estimates for the operators $U_{\epsilon, \lambda}$. \findem

\mb After these remarks, we are going to prove Theorem \ref{thm:general} for our model operators.

\section{The study of $T_{\rho,\lambda,L}$ with Littlewood-Paley square functions.}
\label{section:littlewood}

In this section, we obtain the uniform bounds of Theorem \ref{thm:general} with the Hardy spaces when all the exponents $p_i$ are finite. As we will see in Section \ref{section:carleson}, our model operators can be considered as multilinear Calder\'on-Zygmund operators. Consequently, with similar arguments to those of L. Grafakos and N. Kalton used in \cite{GK2}, we can have boundedness of our operators on the sets of {\em atoms} associated to the considered Hardy spaces.
For several years, many papers (see for example \cite{Bownik,MSV,meyer}) emphasize the following problem : how can we extend a linear operator bounded on the set of atoms to the whole Hardy space ? This abstract question is a really problem and does not admit a general positive answer. For example there is a counter-example in \cite{meyer} for the classical Hardy space.

\gb For this reason, we prefer to describe an other proof, which does not use the atomic decomposition of Hardy spaces. That is why, we are going to directly work with the Littlewood-Paley square functions.

\gb For convenience, we deal only with the bilinear case : $n=3$. First remember the definition of our operator~:
we choose two smooth functions $\Phi^{1}$ and $\Phi^2$ with bounded spectrum and we choose a smooth function $\Psi$ whose the spectrum is included in a corona around $0$. Then we construct the operator
\begin{align*}
\lefteqn{T_{\rho,\lambda,L}(f,g)(x):=} & & \\
 & &  \sum_{k\in\Z} L(k) \int_{\R^d}  \Psi_{2^k}(y) [\Phi^{1}_{\lambda_1 2^k} \ast
f](x-\rho_1\lambda_1 y) [\Phi^{2}_{\lambda_2 2^k}\ast g](x-\rho_2\lambda_{2}y)dy.
\end{align*}

\mb To study this last one, we decompose the two functions $f$ and $g$ with the classical wavelets decomposition~:

\begin{lem} Let $\psi$ a smooth function such that
$$ c(\psi):=\int_0^\infty \left| \widehat{\psi}(t\xi)\right|^{2} \frac{dt}{t}$$
be a nonnegative constant independent with respect to $\xi$ (for example, we can just assume that the function $\psi$ is odd and radial). Then we have the decomposition ~:
\be{lem} f=c(\psi)^{-1}\int_{\R^d} \int_{0}^\infty \langle f, \psi_{t,q}\rangle \psi_{t,q} \frac{dtdq}{t}.\ee
In addition, the integral is absolutely convergent for a function $f\in \s(\R^d)$. 
\end{lem}

\dem The result is well-known for $f\in L^2(\R^{d})$, it is shown in the book \cite{Gra} at the chapter 5.6. When $f\in \s(\R^d)$, integrations by parts give us fast decay for $\langle f, \psi_{t,q}\rangle$ and so permit us to prove the absolute convergence.
\findem

\mb From now, we will choose a smooth function $\psi$ which verifies the assumption of the previous Lemma and whose the spectrum is included in a corona around $0$. We decompose also the two functions $f$ and $g$ with the previous lemma and we have also to study the following quantity~:
\begin{align*} 
\lefteqn{F(k,u,v,q,s,x):=} & & \\
 & &  L(k) \int_{\R^d} \Psi_{2^k}(y) [\Phi^1_{\lambda_1 2^k} \ast \psi_{u,q}](x-\rho_1\lambda_1y) [\Phi^2_{\lambda_2 2^k}
\ast \psi_{v,s}](x-\rho_2\lambda_2 y)dy.
\end{align*} 
With the inverse Fourier transform, we get~:
\begin{align*}
\lefteqn{F(k,u,v,q,s,x) = } & & \\
& & L(k) \int_{\R^{2d}} \widehat{\Psi_{2^p}}(\xi)\widehat{\Phi^1_{\lambda_1
2^k}}((\xi-\alpha)\rho_1^{-1}\lambda_1^{-1})\widehat{\psi_{u,q}}((\xi-\alpha)\rho_1^{-1}\lambda_1^{-1})
e^{i(\xi-\alpha)x\rho_1^{-1}\lambda_1^{-1}}  \\
& &   \widehat{\Phi^2_{\lambda_22^k}}(\lambda_2^{-1}\rho_2^{-1}\alpha) \widehat{\psi_{v,s}}(\rho_2^{-1}\lambda_2^{-1}\alpha)e^{i\alpha x(\rho_2\lambda_2)^{-1}} (\rho_1\rho_2|\lambda_1\lambda_2|)^{-d} d\alpha d\xi.
\end{align*} 
Due to the spectral conditions, we have a dependence for the three frequency parameters~: 
$$F(k,u,v,q,s,x)\neq 0 \Longrightarrow 
\left\{ \begin{array}{l}
\max\{|\rho_1 \lambda_1|u^{-1},|\rho_2\lambda_2|v^{-1}\} \simeq 2^{-k} \\ \textrm{ or } \\
 |\rho_1\lambda_1|^{-1}u\simeq |\rho_2\lambda_2|^{-1}v\leq 2^k \end{array} \right. .$$
In addition the product $\widehat{\Phi^1_{\lambda_1
2^k}}\widehat{\Psi_{u,q}}$ is non vanishing only if $|\lambda_1|u^{-1} \lesssim 2^{-k}$ and similarly for $v$.
As the coefficients $\rho_i$ are bounded by $1$, we are always in the first case i.e.
$$ \max\{|\rho_1 \lambda_1|u^{-1},|\rho_2\lambda_2|v^{-1}\} \simeq 2^{-k}.$$
In addition, we have shown the stronger condition
$$ \max\{|\lambda_1|u^{-1},|\lambda_2|v^{-1}\} \simeq 2^{-k}$$
uniformly with respect to the parameters $\rho_i\in]0,1]$. 

\mb After having study the frequency properties of $F(k,u,v,q,s,x)$, we will remember spatial estimates~:

\begin{prop} \label{coeffs} We have the following estimate~:
$$ F_{k,v,s}(x):=\left|\Phi_{\lambda_2 2^k} \ast \psi_{v,s}(x) \right| \lesssim \frac{\inf\{2^k|\lambda_2|,v\}^d}{v^{2d}}\left(1+\frac{|x-s|}{\max\{|\lambda_2|2^k,v\}}\right)^{-M},$$
 for any exponent $M$ as large as we want. This estimate is uniform with respect to
$k$ and $\lambda_2$. \end{prop}

\dem The proof is essentially written in Appendix K-2 of \cite{Gra} and we only give the sketch of the proof. \\
Let $\tilde{\psi}$ be an other smooth function, whose the spectrum is included in a corona around $0$ and such that
$\widehat{\tilde{\psi}}=1$ on the spectrum of $\psi$. We set $\zeta=\Phi_{\lambda_2 2^k} \ast
\tilde{\psi}_{v,s}$. It is also easy to check that
$$\left|\zeta(x)\right|\lesssim \frac{1}{v^d} \left( 1+ \frac{|x|}{\max\{2^k|\lambda_2|,v\}}\right)^{-M},$$
for all exponent $M$ as large as we want. Due to the spectral properties of $\psi$ and $\tilde\psi$, we get~:
$$F_{k,v,s}(x)= \left|\zeta \ast \psi_{v,s}(x) \right|.$$
Then we can directly estimate the convolution product and prove what we want.\findem

\mb After this study, we decompose our operator~:
\begin{align*}
\lefteqn{T_{\rho,\lambda,L}(f,g)(x)=} & & \\
 & & \sum_{k} \iint_{\R^{2d}}\int_0^\infty \int_0^\infty \langle f, \psi_{u,q}\rangle \langle g, \psi_{v,s}\rangle
 F(k,u,v,q,s,x) \frac{dv}{v}\frac{du}{u} dqds.
\end{align*}

\mb We have seen by a spectral analysis that we can restrict this double integral over $u$ and $v$ on the set 
$$ \left\{ (u,v), \  \max\{|\lambda_1|u^{-1},|\lambda_2|v^{-1}\} \simeq 2^{-k}\right\}.$$ \\
In the study of paraproducts (see paragraph 8.4 of \cite{Gra}), we decompose the product as $fg=\Pi_f(g)+\Pi_g(f)+D(f,g)$ where the two paraproducts and the diagonal terms have different estimates. For the same reasons here we have to singly study the two following terms~:
$T_{\rho,\lambda,L}=T_{\rho,\lambda,L}^{1}+T_{\rho,\lambda,L}^{2}$ with
\begin{align*}
\lefteqn{T_{\rho,\lambda,L}^{1}(f,g)(x):=} & & \\
 & & \sum_{k} \iint_{\R^{2d}}\int_0^\infty \int_0^\infty \langle f, \psi_{u,q}\rangle \langle g, \psi_{v,s}\rangle
 F(k,u,v,q,s,x) {\bf 1}_{A_k}(u,v) \frac{dvdu}{vu} dqds
\end{align*}
and
\begin{align*}
\lefteqn{T_{\rho,\lambda,L}^{2}(f,g)(x):=} & & \\
 & &  \sum_{k} \iint_{\R^{2d}}\int_0^\infty \int_0^\infty \langle f, \psi_{u,q}\rangle \langle g, \psi_{v,s}\rangle
 F(k,u,v,q,s,x) {\bf 1}_{B_k}(u,v) \frac{dvdu}{vu} dqds.
\end{align*}
We write $A_k$ and $B_k$ for the two sets~: \begin{align*}
 \lefteqn{ A_k:=} & & \\
 & & \left\{(u,v), \
\max\left\{|\lambda_1|u^{-1},|\lambda_2|v^{-1}\right\} \simeq 2^{-k} \textrm{ and } \min\{u,v\} \leq C^{-1} \max\{u,v\} \right\} \end{align*}
and
 $$ B_k:=\left\{(u,v), \
\max\left\{|\lambda_1|u^{-1},|\lambda_2|v^{-1}\right\} \simeq 2^{-k} \textrm{ and } u \simeq v \right\}, $$ 
where $C$ is a numerical constant, we later choose. Due to this constant, we can use spectral separation to study $T_{\rho,\lambda,L}^{1}$ with 
 the Littlewood-Paley square functions for $f$ and $g$.

\begin{thm} \label{th4} For $0<\rho_i\leq 1$ and $\lambda_i\in \R^*$, the operator $T_{\rho,\lambda,L}^{1}$ can been continuously extended from
$H^{r_1} \times H^{r_2}$ to $H^{r_3}$, if the exponents $0<r_1,r_2,r_3<\infty$ satisfy the homogeneity condition
$$\frac{1}{r_1} + \frac{1}{r_2} =\frac{1}{r_3}.$$ In addition we control the continuity bounds, uniformly with respect to $\lambda$ and $\rho$ by the quantity 
$$c_N(\Psi)c_N(\Phi^1)c_N(\Phi^2)\|L\|_\infty,$$ 
for $N$ a large enough integer.
\end{thm}

\dem To estimate $T_{\rho,\lambda,L}^{1}$ in the Hardy space $H^{r_3}$, we have to study its square function : $S_\Psi(T_{\rho,\lambda,L}^{1}(f,g))$. We can compute the Fourier transform and get~: 
\begin{align}
\widehat{F(k,u,v,q,s,.)}(\xi)= & L(k) \int_{\R^d} \widehat{\Psi_{2^k}}
(\xi\rho_1\lambda_1-(\rho_2\lambda_2-\rho_1\lambda_1)\alpha)  \nonumber \\
 & \widehat{\psi_{u,q}}(\xi-\alpha)\widehat{\Phi^{1}_{\lambda_1 2^k}}(\xi-\alpha)
\widehat{\Phi^{2}_{\lambda_2 2^k}}(\alpha)\widehat{\psi_{v,s}}(\alpha)d\alpha. \label{fourierF} 
\end{align}
Consequently by writing $\xi=(\xi-\alpha)+\alpha$, the spectrum of $F(k,u,v,q,s,.)$ is contained in 
\be{sommealg} \frac{1}{u} \textrm{spectrum}(\psi)+ \frac{1}{v} \textrm{spectrum}(\psi) \subset \left\{ \xi,\ |\xi|\simeq \max\{u^{-1},v^{-1}\} \right\}. \ee For the last inclusion, we have used a large enough constant $C$ in the definition of the set $A_k$ and so
$$\min\{u,v\} \ll \max\{u,v\}.$$
By symmetry we may assume $u\leq v$ and then by choosing a continuous square function $S_\Psi$, we have
\begin{align*}
S_\Psi(T_{\rho,\lambda,L}^{1}(f,g))(x) \lesssim & \left( \int_{0}^\infty \left[\sum_{k\in\Z} \iint_{\R^{2d}} \int_0^\infty \left|\langle f, \psi_{u,q}\rangle
 \langle g, \psi_{v,s}\rangle\right| \right. \right. \\
 &   \left|\Psi_{u} \ast F(k,u,v,q,s,.) (x)\right| {\bf 1}_{A_k'}(u,v) \frac{dv}{v}
 dqds\Bigg]^2 \frac{du}{u} \Bigg)^{1/2},
 \end{align*}
where $A_k':=\{(u,v)\in A_k, \quad u\ll v\}$. We have also to estimate the following quantity
\begin{align*}
\lefteqn{Q:= \left\| \left(\int_{0}^\infty  \left| \sum_{k}\iint_{\R^{2d}} \int_{0}^\infty  |\langle f, \psi_{u,q}\rangle \langle g, \psi_{v,s}\rangle| \right. \right. \right.} & & \\
 & & \left.  \left|\Psi_{u}\ast F(k,u,v,q,s,.) (x)\right| {\bf 1}_{A_k'}(u,v) \frac{dv}{v} dqds\Bigg|^2 \frac{du}{u} \Bigg)^{1/2}
\right\|_{r_3,dx}.
\end{align*}
With the notations of Proposition \ref{coeffs}, we have that for all $z\in\R^d$
$$ \left|F(k,u,v,q,s,z)\right|\leq |L(k)| \int_{\R^d} \left|\Psi_{2^k} (y) \right| F_{k,u,q}(z-\lambda_1\rho_1y)F_{k,v,s}(z-\rho_2\lambda_2 y) dy.$$
By using the estimates of Proposition \ref{coeffs} and the fast decay of $\Psi$, we get~:
\begin{align*}
\left|F(k,u,v,q,s,z)\right|  & \lesssim   \|L\|_\infty \int_{\R^d} 2^{-dk}\left(1+ \frac{|y|}{
2^k}\right)^{-M} \left(1+ \frac{|z-\rho_1\lambda_1y-q|}{u}\right)^{-M}  \nonumber \\
 &  \quad |\lambda_1|^d 2^{k}u^{-2d} \frac{|\lambda_2|^d 2^{k}}{v^{2d}}\left(1+\frac{|z-\rho_2 \lambda_2 y -s|}{v}\right)^{-M} dy.
\end{align*}
As the parameters $\rho_i \leq 1$ and $(u,v)\in A_k'$, we obtain with an other exponent~:
\begin{align}
\lefteqn{\left|F(k,u,v,q,s,z)\right| \lesssim} & & \\
 & & \|L\|_\infty  \left(1+\frac{|z-s|}{v}\right)^{-M} \left(1+
 \frac{|z-q|}{u}\right)^{-M}\left(\frac{|\lambda_1||\lambda_2|2^{2k}}{u^2v^2}\right)^d. \label{FF}
\end{align}
The exponent $M$ is not always the same, but it always corresponds to an integer as large as we want. Now by estimating the convolution product, we get
\begin{align*}
\lefteqn{\left|\Psi_{u}\ast F(k,u,v,q,s,.) (x)\right| \lesssim} & & \\ 
  & & \|L\|_\infty  \left(1+\frac{|x-s|}{v}\right)^{-M} \left(1+
  \frac{|x-q|}{u}\right)^{-M}\left(\frac{|\lambda_1||\lambda_2|2^{2k}}{u^2v^2}\right)^d. \end{align*}
\mb Computing this estimate in the expression of $Q$, we have
\begin{multline*}
 Q\lesssim \|L\|_\infty \left\| \left\| \sum_{k\in\Z} \iint_{\R^{2d}}\int_0^\infty |\langle f,\psi_{u,q}\rangle\langle g,
\psi_{v,s}\rangle| \left(1+\frac{|x-s|}{v}\right)^{-M} \right.\right. \\
\left.\left. \left(1+\frac{|x-q|}{u}\right)^{-M} {\bf 1}_{A_k'}(u,v) \left(\frac{|\lambda_1||\lambda_2|2^{2k}}{u^2v^2}\right)^d \frac{dv}{v} dqds \right\|_{2,du/u}
\right\|_{r_3,dx}.\end{multline*} 
We change the two variables
$$\frac{x-q}{u} \rightarrow a \quad \textrm{and} \quad \frac{x-s}{v}\rightarrow b,$$
to get
\begin{multline*}
 Q\lesssim \left\| \left\|  \sum_{k} \iint_{\R^{2d}}\int_0^\infty |\langle f,\psi_{u,x-ua}\rangle\langle g,
\psi_{v,x-vb}\rangle|
 \left(1+\left|a\right|\right)^{-M} \left(1+ |b|\right)^{-M}  \right.\right. \\
\left.\left. {\bf 1}_{A_k'}(u,v) \left(\frac{|\lambda_1||\lambda_2|2^{2k}}{uv}\right)^d \frac{dv}{v} dadb \right\|_{2,du/u} \right\|_{r_3,dx}.
\end{multline*}
We write $\psi^{z}$ for $\psi^z(y):= \psi(z-y)$  and $\psi^{z}_t:=(\psi^z)_t$. With these notations, we have~: \begin{multline*} Q\lesssim \left\|
\left\| \sum_{k} \iint_{\R^{2d}}\int_0^\infty |\psi^a_{u} \ast f(x)||\psi^{b}_{v} \ast g(x)| \left(1+\left|a\right|\right)^{-M} \left(1+
|b|\right)^{-M}  \right.\right. \\
\left.\left. {\bf 1}_{A_k'}(u,v)\left(\frac{|\lambda_1||\lambda_2|2^{2k}}{uv}\right)^d \frac{dv}{v} dadb \right\|_{2,du/u} \right\|_{r_3,dx}.
\end{multline*}
The definition of the set $A_k'$ allows us to have a finite number of choice for $k$. Therefore we have
\begin{eqnarray*}
\lefteqn{Q\lesssim \left\| \left\| \iint_{\R^{2d}}\int_0^\infty |\psi^a_{u} \ast f(x)||\psi^{b}_{v} \ast g(x)| \left(1+\left|a\right|\right)^{-M} \left(1+ |b|\right)^{-M} \right.\right.} & &  \\
& & \left.\left.  {\bf 1}_{u\leq C^{-1}v} \left(\frac{|\lambda_1||\lambda_2| \inf\{ |\lambda_1|^{-1}u,|\lambda_2|^{-1}v\}^{2}}{uv}\right)^d \frac{dv}{v} da db
\right\|_{2,du/u} \right\|_{r_3,dx}.
\end{eqnarray*}
Then we use the Cauchy-Schwarz inequality with 
\begin{align*}
(1) & := \int_{Cu}^{\infty} \left(\frac{|\lambda_1||\lambda_2| \inf\{ |\lambda_1|^{-1}u,|\lambda_2|^{-1}v\}^2}{uv}\right)^{2d} \frac{dv}{v} \\
 & \lesssim \int_{Cu}^{\infty} \left(\frac{\inf\{ |\lambda_1|^{-1}u,|\lambda_2|^{-1}v\}}{\max\{ |\lambda_1|^{-1}u,|\lambda_2|^{-1}v\}}\right)^{2d} \frac{dv}{v} \\
 & \lesssim \int_{|\lambda_2\lambda_1^{-1}|u}^{\infty} \left(\frac{|\lambda_1|^{-1}u}{|\lambda_2|^{-1}v} \right)^{2d} \frac{dv}{v} + \int_{Cu}^{\max\{C,|\lambda_2\lambda_1^{-1}|\}u} \left(\frac{|\lambda_2|^{-1}v}{|\lambda_1|^{-1}u} \right)^{2d}
\frac{dv}{v}  \\
& \lesssim 1 +  \left(\frac{|\lambda_2|^{-1}}{|\lambda_1|^{-1}}
\right)^{2d}\left[\max\{C,|\lambda_2\lambda_1^{-1}|\}^{2d}-C^{2d}\right] \\
& \lesssim 1 + \left(\frac{|\lambda_1|}{|\lambda_2|}
\right)^{2d}\left[|\lambda_2\lambda_1^{-1}|^{2d} - \min\{C,|\lambda_2\lambda_1^{-1}|\}^{2d}\right] \\
& \lesssim 1 + \left[1 - \min\{C|\lambda_2^{-1}\lambda_1|,1\}^{2d}\right] \lesssim 1,
\end{align*}
to finally get
$$ Q\lesssim \left\|\left\|  \iint_{\R^{2d}} \left|\psi^{a}_{u} \ast f(x)\right| S_{\psi^{b}}(g)(x)
 \left(1+\left|a\right|\right)^{-M}
 \left(1+ |b|\right)^{-M}dadb\right\|_{2,du/u} \right\|_{r_3,dx}.$$
With Minkowski inequality, we may write the last inequality as
\be{relation} Q\lesssim \left\|\iint_{\R^{2d}} S_{\psi^{a}}(f)(x) S_{\psi^{b}}(g)(x)  \left(1+\left|a\right|\right)^{-M}
 \left(1+ |b|\right)^{-M}dadb \right\|_{r_3,dx}. \ee

\mb We must have a pointwise estimate on the square functions when $r_3<1$, because of the lack for the triangle inequality in $L^{r_3}$. We also use this lemma~:

\begin{lem} \label{ponctuel} Let $\zeta$ be a smooth function satisfying
\be{egal} \forall \xi\neq 0, \qquad \int_0^\infty \left|\widehat{\zeta}(t\xi)\right|^2 \frac{dt}{t}=1. \ee For $\psi$ an other function satisfying
(\ref{egal}) too, we have the pointwise estimate~: for all $r>0$, there exists a constant $C_r$ and an integer $N$ such that
$$\forall f\in\s(\R^d) \qquad S_{\psi}(f)\leq C_r c_N(\zeta)c_N(\psi) \left(\int \left[M^2_r(\zeta_t \ast f )\right]^2
\frac{dt}{t}\right)^{1/2}.$$ Here $M_{HL}$ is the Hardy-Littlewood maximal operator, $M_{HL}^2$ corresponds to $M_{HL}\circ M_{HL}$ and 
$$M^2_r(g)=\left( {M_{HL}}^2(|g|^r)\right)^{1/r}.$$
\end{lem}

\mb Let assume first this Lemma. By applying it with $\psi=\psi^a$ and $\psi=\psi^b$, we get~: \begin{align*}
\lefteqn{Q\lesssim \left\| \iint_{\R^{2d}} (1+|a|)^N (1+|b|)^N\left(\int_0^\infty \left[M^2_r(\zeta_t \ast f )\right]^2
\frac{dt}{t}\right)^{1/2} \right.} & & \\
 & & \left. \left(\int_0^\infty \left[M^2_r(\zeta_t \ast g )\right]^2 \frac{dt}{t}\right)^{1/2} \left(1+|b|\right)^{-M}
 \left(1+ |a|\right)^{-M}dadb\right\|_{r_3}.\end{align*}
Also by choosing a large enough integer $M$, we have~:
$$ Q\lesssim  \left\|\left(\int_0^\infty \left[M^2_r(\zeta_t \ast f )\right]^2
\frac{dt}{t}\right)^{1/2} \left(\int_0^\infty \left[M^2_r(\zeta_t \ast g )\right]^2 \frac{dt}{t}\right)^{1/2}\right\|_{r_3}.$$
With Hölder inequality, we obtain~:
\be{fin} Q\lesssim  \left\|\left(\int_0^\infty \left[M^2_r(\zeta_t \ast f )\right]^2 \frac{dt}{t}\right)^{1/2}\right\|_{r_1}
\left\|\left(\int_0^\infty \left[M^2_r(\zeta_t \ast g )\right]^2 \frac{dt}{t}\right)^{1/2}\right\|_{r_2}. \ee We study only the first term with $r_1$, the other one is identical. By definition,
$$\left\| \left(\int_0^\infty \left[M^2_r(\zeta_t \ast f )\right]^2
\frac{dt}{t}\right)^{1/2} \right\|_{r_1} = \left\|\left(\int_0^\infty \left[M_{HL}^2(|\zeta_t \ast f|^r )\right]^{2/r}
\frac{dt}{t}\right)^{r/2}\right\|_{r_1/r}^{1/r}.$$ For $r$ small enough such that $\min\{r_1/r,2/r\}> 1$, the Fefferman-Stein inequality (Theorem 4.6.6 of \cite{Gra}) in $L^{2/r}$ applied to the operator $M_{HL}^2$ gives us~:
$$\left\| \left(\int_0^\infty \left[M^2_r(\zeta_t \ast f )\right]^2
\frac{dt}{t}\right)^{1/2} \right\|_{r_1} \lesssim \left\|\left(\int_0^\infty \left[|\zeta_t \ast f|^r \right]^{2/r}
\frac{dt}{t}\right)^{r/2}\right\|_{r_1/r}^{1/r}.$$ In other words~:
$$\left\| \left(\int_0^\infty \left[M^2_r(\zeta_t \ast f )\right]^2
\frac{dt}{t}\right)^{1/2} \right\|_{r_1} \lesssim \left\|S_{\zeta}(f) \right\|_{r_1}.$$ By replacing this estimate in (\ref{fin}), we obtain the desired result~:
$$Q\lesssim \|f\|_{H^{r_1}} \|g\|_{H^{r_2}}$$
uniformly with respect to $\lambda_i\neq 0$ and $0<\rho_i\leq 1$.
\findem

\mb We have now to show the Lemma \ref{ponctuel}. This Lemma is ``quite classical'', it permits to understand for example that Definition of Hardy spaces (Definition \ref{hardyspace}) does not depend on the used function $\Psi$. It is almost proved in a discrete version in the book \cite{Gra}, from which we take the notations. We only give the sketch of the proof.

\dem We define the maximal operator~:
$$M_{b,t}(f,\psi)=\sup_{y\in\R^d} \frac{\left|\psi_t\ast f(x-y)\right| }{\left(1+t^{-1}|y|\right)^b}.$$
Then it is obvious that \
\be{eq1lem} \left|\psi_t \ast f\right|\leq M_{b,t}(f,\psi).\ee In the proof of Theorem 6.5.6. of \cite{Gra}, one may choose a function $\Theta$ satisfying (\ref{egal}) and such that $\hat{\Theta}\geq
0$. Then it is shown that ~:
$$M_{b,t}(f,\psi)\lesssim \int  c_N(\psi)c_N(\Theta) \inf\{|t-s|,|t-s|^{-1}\} M_{b,s}(f,\Theta)\frac{ds}{s}.$$
Consequently with (\ref{eq1lem}), we get~:
\be{equalem} S_{\psi}(f)\lesssim c_N(\psi)c_N(\Theta) \left(\int \left[M_{b,s}(f,\Theta)\right]^2\frac{ds}{s}\right)^{1/2}.\ee Now Lemma 6.5.3. of \cite{Gra} with $b=n/r$ gives us, 
$$M_{b,t}(f,\Theta)\lesssim M_r(\Theta_t \ast f).$$
To substitute the function $\zeta$ to the function $\Theta$, we use the spectral condition and the fact that
$$\Theta_t \ast f=\int_{1/2}^2 \Theta_t\ast \zeta_{tu} \ast f \frac{du}{u}.$$
Then with the estimate (6.5.8) of \cite{Gra}~:
$$\left|\Theta_t\ast \zeta_{tu} \ast f(x)\right|\lesssim c_N(\Theta) M_r(\zeta_t \ast f) (x),$$
we obtain that
$$M_{b,t}(f,\Theta)\lesssim M_r(\Theta_t \ast f) \lesssim c_N(\Theta)c_N(\zeta) M_r M_r(\zeta_t \ast f) \lesssim c_N(\zeta) M^2_r(\zeta_t \ast
f).$$ By computing this estimate in (\ref{equalem}), we get the Lemma. \findem

\gb To finish the study of $T_{\rho,\lambda,L}$, we have to estimate the second operator
$T_{\rho,\lambda,L}^{2}$~:

\begin{thm} \label{th6} For $0<\rho_i\leq 1$ and $\lambda_i\neq 0$, the operator $T_{\rho,\lambda,L}^{2}$ is continuous from $H^{r_1} \times H^{r_2}$ to $L^{r_3}$, for all exponents
$0<r_1,r_2,r_3<\infty$ satisfying
$$\frac{1}{r_1} + \frac{1}{r_2} =\frac{1}{r_3}.$$ 
In addition we may control the continuity bounds, uniformly with respect to $\lambda$ and $\rho$ by the quantity 
$$c_N(\Psi)c_N(\psi)c_N(\Phi^1)c_N(\Phi^2)\|L\|_\infty,$$ 
for $N$ a large enough integer.
\end{thm}

\dem The operator $T_{\rho,\lambda,L}^{2}$ is defined as
\begin{align*}  
\lefteqn{T_{\rho,\lambda,L}^{2}(f,g):=} & & \\
 & &  \sum_{k\in\Z} \iint_{\R^{2d}} \int_0^\infty \int_0^\infty  \langle f, \psi_{u,q}\rangle \langle g, \psi_{v,s}\rangle
 F(k,u,v,q,s,x) {\bf 1}_{B_k}(u,v) \frac{dvdu}{vu} dqds,
\end{align*}
with $$B_k:=\left\{(u,v), \ \max\{|\lambda_1|u^{-1},|\lambda_2|v^{-1}\} \simeq 2^{-k} \textrm{ and } u \simeq v \right\}.$$ By using the previous estimate, we have to control
\begin{align*}
Q & :=\left\|T_{\rho,\lambda,L}^{2}(f,g)\right\|_{r_3} \\
 & \lesssim \left\|\sum_{k\in\Z} \iint_{\R^{2d}} \int_0^\infty \int_0^\infty |\langle
f,\psi_{u,q}\rangle\langle g, \psi_{v,s}\rangle| \left(1+ \frac{|x-s|}{v}\right)^{-M} \right. \\
 &  \qquad \qquad \left.  \left(1+
\frac{|x-q|}{u}\right)^{-M}{\bf 1}_{B_k}(u,v) \left(\frac{|\lambda_1||\lambda_2|2^{2k}}{u^2v^2}\right)^d  \frac{dvdu}{uv} dqds
\right\|_{r_3,dx}.
\end{align*}
In this case (\ref{sommealg}) is not satisfied. We compute the same changes of variables as in the end of the proof for Theorem \ref{th4}. By using Cauchy-Schwarz inequality and the definition of the set $B_k$, we obtain the same estimate as (\ref{relation}) and so we can conclude by the same arguments as before.
\findem

\gb Finally we get the following result ~:
\begin{thm}  \label{thm:final1} Let $0<\rho_i\leq 1$ and $\lambda_i\neq 0$ be reals, then the operator $T_{\rho,\lambda,L}$ is continuous from $H^{r_1} \times H^{r_2}$ to $L^{r_3}$ for all exponents
$0<r_1,r_2,r_3<\infty$ satisfying 
$$\frac{1}{r_1} + \frac{1}{r_2} =\frac{1}{r_3}.$$ In addition we can estimate the continuity bound uniformly with respect to $\lambda$ and $\rho$ by the quantity 
$$c_N(\Psi)c_N(\psi)c_N(\Phi^1)c_N(\Phi^2)\|L\|_\infty,$$ 
for $N$ a large enough integer.
\end{thm}

\dem We have decomposed the operator $T_{\rho,\lambda,L}$ as
 $$T_{\rho,\lambda,L}=T_{\rho,\lambda,L}^{1}+T_{\rho,\lambda,L}^{2}.$$ The embedding $H^{r_3} \hookrightarrow L^{r_3}$
(see Theorem 2.5 of \cite{Uchi}), Theorems \ref{th4} and \ref{th6} allow us to prove the desired result.
\findem

\mb This result proves the first part of Theorem \ref{thm:general} : under the assumption 1-) we have uniform estimates.

\mb In the next section, we are going to prove a similar result for some infinite exponents with the concept of Carleson measure and ideas based on Calder\'on-Zygmund theory.

\section{The study of $T_{\rho,\lambda,L}$ with Carleson measures and Calder\'on-Zygmund decompositions .}
\label{section:carleson}

We use ideas of the book \cite{cm}, where R. Coifman and Y. Meyer have already studied paraproducts with a Carleson measure. We adapt here their arguments to our model operators. As we have seen in Remark \ref{remsym}, our operators permit us to understand all the ``different kinds'' of paraproducts. In \cite{cm}, the authors studied only one ``kind'' of paraproducts (which with other and extra arguments is sufficient to study the other ones). \\
That is why the use of our model operators seems interesting as we obtain a (only one) direct proof simultaneously for all the paraproducts.

\mb We will (for convenience) work on the continuous version of them~: 
\begin{eqnarray*}
U_{\rho,\lambda,L}(f_1,..,f_{n})(x):= \int_0^\infty L(t) \int_{\R^d}  \Psi_{t}(y) \prod_{i=1}^{n} \left[ \Phi^{i}_{\lambda_i t}\ast f_i \right](x-\rho_i\lambda_iy)dy\frac{dt}{t},
\end{eqnarray*}
where $L$ is a bounded measurable function. By symmetry, we can assume that $\lambda_{n}$ satisfies~: \be{conditionmu}
|\lambda_{n}|= \max\{|\lambda_i|,\ 1\leq i \leq n\}.\ee In this case we have the following result~:

\begin{thm} \label{th7}  Under the assumption (\ref{conditionmu}), the operator $U_{\rho,\lambda,L}$ can be continuously extended from $(L^\infty)^{\otimes n-1} \times L^2$ to $L^2$. In addition the continuity bound is controlled  by $\|L\|_\infty$ with $c_{M}(\Psi)$ and $c_M(\Phi^{i})$ (for a large enough integer $M$) independently with respect to the parameters $\rho_i\in]0,1]$ and $\lambda$. 
\end{thm}

\dem  By symmetry on the $(n-1)$ first coordinates, we can assume that
\be{conditionmu2} |\lambda_1| := \min \{|\lambda_j|, 1\leq j \leq n-1\}. \ee
In $\cite{cm}$ (Chap. VI prop 3), the following result (that we call the $(*)$-result) is shown~: the operator $V$ is continuous from $(L^\infty)^{\otimes n-1} \times L^{2}$ to $L^2$, where $V$ is defined by
\begin{align*}
V(f_1,..,f_{n})(x) & := U_{(\rho_1,0,..,0),\lambda,L}(f_1,..,f_{n})(x)\\
 & = \int_0^\infty  \left[ \Psi_{\rho_1\lambda_1t}\ast \Phi^{1}_{\lambda_1t} \ast f_1\right] (x)
\prod_{j=2}^{n} \Phi^{j}_{\lambda_j t}\ast f_j(x) \frac{L(t)}{t} dt.\end{align*} The estimate on $V$ is independent on
$\lambda$ and $\rho_1$ due to the assumptions (\ref{conditionmu}) and (\ref{conditionmu2}). Our idea is also to disturb
our coefficients $(\rho_j)_{2\leq j\leq n}$ and to bring them to $0$.  We temporarily forget in the notation the dependence on $\rho,\lambda$ and $L$, by writing~:
$$U_{\rho,\lambda,L}(f_1,..,f_{n})=V(f_1,..,f_{n})+ \sum_{\genfrac{}{}{0pt}{}{J\subset\{2,..,n\} }{J\neq \varnothing} } \int_{\genfrac{}{}{0pt}{}{0\leq s_j \leq \rho_j}{ j\in J}} W_{s,J}(f_1,..f_{n}) ds,$$
where
\begin{align*} \lefteqn{W_{s,J}(f_1,..,f_{n})(x):= \int_0^\infty \int_{\R^d}  \Psi_{t}(y) \Phi^{1}_{\lambda_1 t}\ast f_1(x+\rho_1\lambda_1 y)} & &  \\
 & & \prod_{\genfrac{}{}{0pt}{}{j=2}{j\in J}}^{n} \frac{\lambda_j y}{\lambda_j t}.( \nabla \Phi^{j})_{\lambda_j t}\ast f_j(x+s_j\lambda_j y) \Pi_J(x,y,t) \frac{L(t)}{t} dydt
\end{align*}
with
$$\Pi_J(x,y,t):= \prod_{\genfrac{}{}{0pt}{}{j=2}{j\in J^{c}}}^{n} (\Phi^{j})_{\lambda_j t}\ast f_j(x+\lambda_j y).$$

\mb Since $V$ is estimated by the $(*)$-result, $\rho_i\leq 1$ and the set of $J$ being finite, we have only to bound the operators $W_{s,J}$. We now decompose the gradient in the $d$ coordinates~:
\begin{align*}
\lefteqn{W_{s,J}(f_1,..,f_{n})(x)=\sum_{l\in\{1,..,d\}^{|J|}} \int_0^\infty \int_{\R^d}  \Psi_{t}(y) \Phi^{1}_{\lambda_1 t}\ast f_1(x+\rho_1\lambda_1 y) } & & \\
 & & \prod_{\genfrac{}{}{0pt}{}{j=2}{j\in J}}^{n}\frac{y_{l_j}}{t} \left[ ( \partial_{x_{l_j}} \Phi^{j})_{\lambda_j t}\ast f_j \right] (x+s_j\lambda_j y) \Pi_J(x,y,t) \frac{L(t)}{t} dydt.
  \end{align*}
By setting
 $$\Theta^{l}(x)=\left[\prod_{j\in J} x_{l_j} \right] \Psi(x),$$
the function $\Theta^{l}$ is again a smooth function whose the spectrum is far away from $0$. We have
\begin{align*}
\lefteqn{W_{s,J}(f_1,..,f_{n-1})(x)=\sum_{l} \int_0^\infty \int_{\R^d}  \Theta^{l}_{t}(y) \left[ \zeta^{1}_{\lambda_1 t}\ast f_1 \right](x+\rho_1\lambda_1 y) } & & \\
 & &  \prod_{j\in J} \left[ (\partial_{x_{l_j}}\Phi^{j})_{\lambda_j t}\ast
 f_j\right] (x+\lambda_j s_j y) \Pi_J(x,y,t) \frac{L(t)}{t} dydt.
\end{align*}

\mb Now the interest of this operation is that $J$ being not empty there exists an index $j$ for which $\widehat{\partial_{x_{l_j}}\Phi^{j}}(0)=0$ (what is false for the initial function $\Phi^{j}$)\footnote{\label{rappel2} In this case, we know from \cite{cm} that $\|(\partial_{x_{l_j}}\Phi^{j})_t\ast
f\|_{2,\frac{dxdt}{t}}\lesssim \|f\|_2$}. The $(*)$-result for $V$ is based on the following quadratic estimate (due to the notion of Carleson measure, see \cite{cm})~: 
$$\|\left[\Psi_{\rho_1\lambda_1t}\ast f\right](x) \left[ \zeta^{2}_{\lambda_2t} \ast g\right](x)\|_{2,\frac{dtdx}{t}} \lesssim \|f\|_\infty\|g\|_2,$$
uniformly on $\lambda$ and $\rho$ for $0<\rho_i\leq 1$ and $0<|\lambda_1|\leq |\lambda_2|$. We are going also to produce the same proof for our operator $W_{s,J}$. We have to show a quadratic estimate : for an index
$l\in\{1,..,d\}^{|J|}$~:

\begin{align}
(\spadesuit)_l & :=\left\|\int_{\R^d} \Theta^{l}_{t}(y) \Phi^{1}_{\lambda_1 t} \ast
f_{1}(x-\rho_1\lambda_1y) \right. \nonumber \\
 & \quad \left. \prod_{j\in J} (\partial_{x_{l_j}} \Phi^{j})_{\lambda_j t}\ast f_j(x+\lambda_j s_j y)\prod_{j\in J^c}
\Phi_{\lambda_j t}\ast f_j(x+\lambda_j y) dy \right\|_{2,\frac{dxdt}{t}} \nonumber \\
  & \lesssim \|f_{n}\|_2 \prod_{i=1}^{n-1} \|f_i\|_\infty. \label{amontrer}
\end{align}

\mb $\ast$ First case : $n\in J$. The convolution operators are bounded on $L^\infty$, so we get~:
 $$(\spadesuit)_l \leq \prod_{i\neq n}\|f_i\|_\infty \left\| \int_{\R^d} |\Theta^{l}_{t}(y)| \left|(\partial_{x_{l_{n}}} \Phi^{n})_{\lambda_{n} t}\ast f_{n} (x-\lambda_{n} s_{n} y) \right|\ dy\right\|_{2,\frac{dxdt}{t}}.$$
We use Minkowski inequality for the norm in $L^2(dt/t)$ and after we can compute the integral over $y$. Then with the reminder
(\ref{rappel2}), we get (\ref{amontrer}). \\
$\ast$ Second case : $n\in J^c$. We use a Carleson estimate by keeping an other function $f_k$ with $k\in J$ (due to $J\neq \emptyset$).
\begin{align*}
\lefteqn{(\spadesuit)_l \leq \prod_{i\neq\{n,k\}}\|f_i\|_\infty \left\| \int_{\R^d} |\Theta^{l}_{t}(y)| \left| \Phi^{n}_{\lambda_{n} t} \ast f_{n}(x-\lambda_{n} y)\right| \right.} & &  \\
 & & \left|(\partial_{x_{l_k}} \Phi^{k})_{\lambda_k t} \ast f_k (x-\lambda_k s_k)\right| dy\Bigg\|_{2,\frac{dxdt}{t}}.
\end{align*}
After changing the variable on $y$, we have~:
\begin{align*}
\lefteqn{(\spadesuit)_l \leq \prod_{i\neq\{n,k\}}\|f_i\|_\infty \left\| \int_{\R^d} |\Theta^{l}(y)| \left|\Phi^{n}_{\lambda_{n} t} \ast f_{n}(x-\lambda_{n} ty)
\right| \right. } & & \\
 & &  \left|(\partial_{x_{l_k}} \Phi^{k})_{\lambda_k t} \ast f_k (x-\lambda_k t s_k y)\right|\ dy\bigg\|_{2,\frac{dxdt}{t}}.
\end{align*}
We write $\Phi^{k,a}$ for $\Phi^{k}(\cdot-a)$ and also get
\begin{align*}
\lefteqn{(\spadesuit)_l\leq \prod_{i\neq\{n,k\}}\|f_i\|_\infty} & & \\
 & &  \left\| \int_{\R^d} |\Theta^{l}(y)| \left|\Phi^{n,y}_{\lambda_{n} t} \ast f_{n}(x)\right|
\left|(\partial_{x_{l_k}}\Phi^{k,s_ky})_{\lambda_k t} \ast f_k (x)\right|\ dy\right\|_{2,\frac{dxdt}{t}}.
\end{align*}
 We now use Minkowski inequality on the measure $dxdt/t$ and after the Carleson estimate to finally obtain~:
$$(\spadesuit)_l \leq \prod_{i\neq n}\|f_i\|_\infty \|f_{n}\|_2 \int |\Theta^{l}(y)| (1+|s_ky|)^{2d+2} dy.$$
The function $\Theta^{l}$ is smooth and $0\leq s_k\leq 1$, consequently we have shown $(\spadesuit)_l$ in this last case. All the estimates are uniform on $\lambda$ due to $|\lambda_k|\leq |\lambda_{n}|$.

\mb Hence (\ref{amontrer}) is shown in the two cases. We have now just to copy the proof of the $(*)$-result of \cite{cm} by putting the previous quadratic estimate instead of the Carleson estimate. The details of this part of the proof are left to the reader. \findem

\gb As for ``classical'' Calder\'on-Zygmund operators, we use a Calder\'on-Zygmund decomposition to obtain continuity results with other Lebesgue exponents. Our multilinear operators are multilinear Calder\'on-Zygmund operators (as defined in  \cite{GT}), however the bounds are not uniformly controlled with respect to $\lambda_i$. By using the main result of \cite{GT}, we obtain our desired continuities for $U_{\rho,\lambda,L}$ with a certain dependence on $\lambda$. The rest of this section is based on an improvement of the ``classical'' arguments, adapted to our problem.

\begin{df} \label{associe} A function $K$ defined on $\R^d \times \R^d \setminus \{(x,x),x\in\R^d \}$ is called a ``standard kernel of order $N$'' if for all $x\neq y$
$$\forall \alpha,\beta \in\{1,..,d\},\ |\alpha|,|\beta|\leq N  \qquad \left|\partial^\alpha_{x} \partial^{\beta}_{y} K(x,y)\right| \leq A_{\alpha,\beta} \frac{1}{|x-y|^{d+|\beta|+|\alpha|}}. $$
A linear operator $T$, continuously acting from $\s(\R^d)$ to $\s'(\R^d)$ and satisfying the integral representation
$$ \forall f\in C^\infty_0,\ \forall x \notin \textrm{supp}(f) \qquad T(f)(x) = \int_{\R^d} K(x,y) f(y) dy, $$
is said to be associated to the kernel $K$.
Such an operator is called a "Calder\'on-Zygmund operator of order $N$" if it is bounded on $L^2(\R^d)$ and associated to a "standard kernel of order $N$".
\end{df}

\mb We have also the well-known following proposition (see for example the book \cite{Gra})~:

\begin{prop} \label{propnoyau} Let $T$ be a Calder\'on-Zygmund operator of order $N$. Then $T$ admits a continuous extension from $L^p$ to $L^p$ for $1<p\leq 2$, from $L^1$ to $L^{1,\infty}$ and from $H^{p}$ to $L^{p}$ for $d/(N+d)<p\leq 1$. In addition the continuity bounds only depend on the constants $\|T\|_{L^2\rightarrow L^2}$ and $(A_{0,\beta})_\beta$.
\end{prop}

\begin{rem} The other constants $A_{\alpha,\beta}$ with $\alpha\neq 0$ are useful to study the dual operator $T^{*}$ and also to get boundedness for $T$ on $L^{q}$ with $2<q<\infty$.
\end{rem}

\mb We will use this proposition for our problem.

\begin{prop} \label{para2}  Let $f_1,..,f_{n-1}$ be smooth fixed functions (considered in $L^\infty$). Then the operator~:
$$V := f_{n} \rightarrow U_{\rho,\lambda,L}(f_1,..,f_{n})$$
is a Calder\'on-Zygmund operator at any order. In addition the constants $A(0,\beta)$ are uniformly bounded with respect to $\lambda$ and $\rho$ for $0<\rho_i\leq 1$.
\end{prop}

\dem The boundedness on $L^2$ of $V$ is given by Theorem \ref{th7}. We have only to check the desired estimates on the kernel. Let $K$ be the kernel of $V$, which is given by~:
\begin{align*}
 \lefteqn{K(x,z) =} & & \\
 & & \int_0^\infty \int_{\R^d}  \Psi_{t}(y) \prod_{j=1}^{n-1} \Phi^{j}_{\lambda_{j}t} \ast f_j (x-\rho_{j}\lambda_{j}y)
  \Phi^{n}_{\lambda_{n}t}(x-\rho_{n}\lambda_{n}y-z)  \frac{L(t)}{t} dydt.
\end{align*}
We can differentiate the kernel and directly obtain~:
\begin{eqnarray*}
\lefteqn{\left| \partial^\alpha_{x} \partial_z^\beta K(x,z)\right| \lesssim \frac{1}{(\min_i |\lambda_i| )^{|\alpha|}} \int_0^\infty
\prod_{j=1}^{n-1} \|f_j\|_\infty \|L\|_\infty} & & \\
& &    \int_{\R^d} \left(1+\frac{|y|}{t}\right)^{-M} \left(1+\frac{|x-\rho_{n}\lambda_{n}y-z|}{\lambda_{n}
t}\right)^{-M}|\lambda_{n}t|^{-d-|\beta|}  \frac{dydt}{t^{|\alpha|+d+1}}.
\end{eqnarray*}
By using $\rho_{n-1}\leq 1$, we get~:
\begin{align*}
 \lefteqn{\left|\partial_x^\alpha \partial_z^\beta K(x,z)\right| \lesssim} & & \\ 
 & & \left(\frac{|\lambda_{n}|}{\min_i |\lambda_i|}\right)^{|\alpha|} \int_0^\infty \left(\frac{1}{|\lambda_{n}|t }\right)^{|\beta|+|\alpha|}
\left(1+\frac{|x-z|}{\lambda_{n} t}\right)^{-M}|\lambda_{n}t|^{-d} \|L\|_\infty  \frac{dt}{t}.
\end{align*}
For $M$ a large enough exponent, we can conclude that
$$ \left|\partial_x^\alpha \partial_z^\beta K(x,z)\right| \lesssim \left(\frac{|\lambda_{n}|}{\min_i |\lambda_i| }\right)^{|\alpha|}
  |x-z|^{-d-|\alpha|-|\beta|} \|L\|_\infty.$$
We have also the desired estimates on the kernel and for $\alpha=0$ the estimates do not depend on $\lambda$. \findem

\mb With the two previous propositions, we get the following corollary.

\begin{cor} \label{cor1} The operator $U_{\rho,\lambda,L}$ is continuous from $(L^\infty)^{\otimes(n-1)} \times H^{p}$ into $L^{p}$ for all exponent $0<p\leq 2$ and from $(L^\infty)^{\otimes(n-1)} \times L^{1}$ into $L^{1,\infty}$. The continuity bounds are uniformly controlled with respect to $0<\rho_i\leq 1$ and $\lambda$ satisfying (\ref{conditionmu}).
\end{cor}

\mb Here we do not know if a similar result for $p>2$ is possible. 

\mb Now we would like to get continuities with finite exponents instead of infinite exponents. To do this, we first prove the abstract following result~:

\begin{thm} \label{L1faible}  Let $T$ be a linear operator, continuously acting from $L^{p_1}$ to $L^{p}$ with $1<p_1\leq \infty$ and $0<p\leq p_1$. We set $r>0$ the exponent defined by
 $$\frac{1}{p}=\frac{1}{p_1}+\frac{1}{r}.$$
We assume that $T$ is associated to a kernel $K$ (see definition \ref{associe}) satisfying
$$ \forall \alpha, |\alpha|\leq 1 \qquad \left|\partial_{z}^\alpha K(x,z) \right| \leq \frac{1}{|x-z|^{d+|\alpha|}} h(x)$$
with a function $h\in L^{r,\infty}$. Then $T$ can be continuously extended from $L^{q_1}$ to $L^{q,\infty}$ for all exponents $(q_1,q)$ such that
$$\frac{1}{q}=\frac{1}{q_1}+\frac{1}{r} \qquad \textrm{and} \qquad  1\leq q_1 \leq p_1.$$
In addition the continuity bounds are controlled by $\|h\|_{r,\infty}$.
By real interpolation, we obtain the strong type $(q_1,q)$ when $1< q_1 \leq p_1$.
\end{thm}

\dem We follow the ``classical proof'' for $r=\infty$ and $h={\bf 1}_{\R^d}$. So let $f$ be a normalized function of $L^{q}$~: $\|f\|_{q}=1$. We want to show
\be{faible} \left|\left\{ x, \ |T(f)|>\alpha \right\} \right| \lesssim \alpha^{-q}.\ee
We use a Calder\'on-Zygmund decomposition of the function $f$ at the scale $\alpha^{q/q_1}$. We have also the following decomposition
$$f=g+b,$$
with a ``good'' function $g$ and a ``bad'' function $b$ satisfying~:
\begin{eqnarray*}
 & \displaystyle \|g\|_{q_1} \lesssim \|f\|_{q_1}=1, \qquad \|g\|_\infty \lesssim \alpha^{q/q_1}, & \\
 & \displaystyle  b=\sum_{k} b_{k}, \qquad \textrm{supp}(b_{k}) \subset Q_{k},  & \\
 & \displaystyle \|b_{k}\|_{q_1} \lesssim \alpha^{q/q_1}|Q_{k}|^{1/q_1}, \qquad \int b_{k} = 0,   &  \\
 & \displaystyle \sum_{k} |Q_{k}| \lesssim \alpha^{-q} \|f\|_{q_1}^{q_1} \lesssim \alpha^{-q} \quad \textrm{and} \quad \sum_{k} {\bf 1}_{10 Q_k} \lesssim 1. &
\end{eqnarray*}
The $(Q_k)_k$ is a collection of balls (of $\R^d$), associated to the ``bad'' function $b$. By linearity, we have
$$ \left|\left\{ x, \ |T(f)|>\alpha \right\} \right| \leq \left|\left\{ x, \ |T(g)|>\alpha/2 \right\} \right| + \left|\left\{ x, \ |T(b)|>\alpha/2 \right\} \right|.$$
$1-)$ The case of the function $g$. \\
This is the easiest case. We use the continuity of $T$ with the exponents $p_1$ and $p$ to get
\begin{align*}
\left|\left\{ x, \ |T(g)|>\alpha/{2} \right\} \right| & \lesssim \alpha^{-p} \left\| T(g)\right\|_{p}^{p} \\
& \lesssim \alpha^{-p} \|g\|_{p_1}^{p}.
\end{align*}
By the assumption on $g$ and the fact that $q_1\leq p_1$~:
$$\|g\|_{p_1} \lesssim \|g\|_{q_1}^{q_1/p_1} \|g\|_\infty^{1-q_1/p_1} \lesssim \alpha^{q(1-q_1/p_1)/q_1}.$$
We also obtain
\begin{align*}
\left|\left\{ x, \ |T(g)|>\alpha/{2} \right\} \right| & \lesssim \alpha^{-p} \alpha^{p q(1-q_1/p_1)/q_1} \\
 & \lesssim \alpha^{-p} \alpha^{p q(1/q_1-1/p_1)} \lesssim \alpha^{-p} \alpha^{pq(1/q-1/p)} \\
 & \lesssim \alpha^{-q},
\end{align*}
which corresponds to the desired result (\ref{faible}). \\
$2-)$ The case of the function $b$. \\
First we have
$$\left| \bigcup_{k} 5 Q_{k} \right| \lesssim \alpha^{-q}.$$
In order to show (\ref{faible}), we can also assume that $x\in \cap_{k} (5 Q_{k})^{c}$ and just estimate
$$\left|\left\{ x\in \cap_{k} (5 Q_{k})^{c}, \ |T(b)|>\alpha/2 \right\} \right|.$$
Let also $x$ be fixed and use
$$ |T(b)(x)| \leq \sum_{k\geq 0} |T(b_{k})(x)|.$$
With the vanishing moment of the function $b_k$, we have~:
$$T(b_{k})(x) = \int K(x,z) b_{k}(z) dz = \int \left[K(x,z)-K(x,c_{k})\right] b_{k}(z) dz.$$
Here we write $c_{k}$ for the center of the cube $Q_{k}$. As $x$ is far away the support of $b_k$, the integral representation has really a sense. Then by using the estimates of the kernel, we have
\begin{align*}
 \left|K(x,z)-K(x,c_{k})\right| & \lesssim |z-c_k| \int_{0}^{1} \left|\nabla K(x,z+t(c_k-z))\right| dt \\  
  & \lesssim |z-c_k| \int_{0}^{1} \frac{1}{|x-c_k|^{d+1}} h(x) dt  \\
  & \lesssim \frac{|z-c_{k}|}{|x-c_{k}|^{d+1}}h(x).
\end{align*}
Therefore
$$\left|T(b_{k})(x)\right| \lesssim h(x) \int |Q_{k}|^{1/d} \frac{1}{|x-c_k|^{d+1}} b_{k}(z) dz.$$
With
$$\|b_k\|_{1} \leq |Q_k|^{1-1/q_1}\|b_{k}\|_{q_1} \lesssim |Q_k|\alpha^{q/q_1},$$
we obtain
$$\left|T(b_{k})(x)\right| \lesssim h(x) |Q_{k}|^{1+1/d} \alpha^{q/q_1} \frac{1}{|x-c_k|^{d+1}}. $$
By computing the sum over the index $k$, we finally have
$$\left|T(b)(x)\right| \lesssim h(x) \alpha^{q/q_1} \sum_{k} \frac{1}{\left(1+\frac{|x-c_k|}{|Q_k|^{1/d}}\right)^{d+1}}. $$
We find the Marcinkiewicz function associated to the collection $(Q_k)$, we write it $M_{(Q_k)_k}$. So we have
$$\left|T(b)(x)\right| \lesssim h(x) \alpha^{q/q_1} M_{(Q_k)_k}(x). $$
However the collection $(10Q_k)_k$ is a bounded covering on the whole space, so we know (see \cite{S}) that for $1\leq q_1 <\infty$, $M_{(Q_k)_k}$ is of weak type $(q_1,q_1)$. By using Hölder inequality on the weak Lebesgue spaces, we get 
$$\left\|T(b)\right\|_{q,\infty} \lesssim \|h\|_{r,\infty} \alpha^{q/q_{1}} \left\|M_{(Q_k)_k}(x)\right\|_{q_1,\infty} \lesssim |\cup Q_{k}|^{1/q_1} \|h\|_{r,\infty} \alpha^{q/q_{1}}.$$
We obtain also the desired estimate~:
$$\left\|T(b)\right\|_{q,\infty} \lesssim \|h\|_{r,\infty}.$$
\findem

\mb We now prove a similar result for the Hardy spaces~:

\begin{thm} \label{hardy}  Let $T$ be a linear operator, continuously acting from $L^{p_1}$ to $L^{p}$ with $1<p_1\leq \infty$ and $p\leq p_1$. We set $r>0$ the exponent satisfying
 $$\frac{1}{p}=\frac{1}{p_1}+\frac{1}{r}.$$
We assume that $T$ is associated (see Definition \ref{associe}) to a kernel $K$ verifying~:
$$ \forall \alpha, |\alpha|\leq N \qquad \left|\partial_{z}^\alpha K(x,z) \right| \leq \frac{1}{|x-z|^{d+|\alpha|}} h(x)$$
with a function $h\in L^{r,\infty}$. Then for all exponents $(q_1,q)$ such that
$$\frac{1}{q}=\frac{1}{q_1}+\frac{1}{r} \qquad \textrm{and} \qquad  \frac{d}{N+d}\leq q_1\leq 1,$$
there is a constant $C$ such that for all  atoms $a\in H^{q_1}$  
 \be{hyphardy} \left\| T(a) \right\|_{L^{q,\infty}} \leq C . \ee
In addition the continuity bounds are controlled by $\|h\|_{r,\infty}$. By real interpolation, $T$ can be continuously extended from $H^{q_1}$ to $L^{q}$ when $\frac{d}{N+d} < q_1\leq 1$.
\end{thm}

\dem We use the atomic decomposition of the Hardy spaces $H^{q_1}$ (See Theorem 6.6.10 of \cite{Gra}). Let $a$ be an atom of $H^{q_1}$, that is meaning there exists a cube $Q$ such that
\begin{align*}
& \textrm{supp}(a)\subset Q, \qquad \|a\|_{2} \leq |Q|^{1/2-1/{q_1}} & \\
& \forall \alpha, |\alpha|\leq [\frac{d}{q_1}-d], \qquad \int x^{\alpha} a(x) dx=0. &
\end{align*}
We write $[]$ for the integer part. We want to estimate $T(a)$. Assume first that $x\in (5Q)^{c}$.
By assumption $q_1>d(N+d)^{-1}$ so $N_{q_1}:=[\frac{d}{q_1}-d]\leq N-1$. We have also
\begin{align*}
T(a)(x) & = \int_{\R^d} K(x,y) a(y) dy \\
 & = \int_{\R^d} \left[K(x,y)- \sum_{|\alpha|\leq N_{q_1}} \frac{(y-c(Q))^\alpha}{\alpha!} \partial^\alpha_{y} K(x,c(Q))\right] a(y) dy,
 \end{align*}
where $c(Q)$ is the center of the cube $Q$. We can estimate the difference between the square brackets by
\begin{align*}
 \left|K(x,y)- \sum_{|\alpha|\leq N_{q_1}} \frac{(y-c(Q))^\alpha}{\alpha!} \partial^\alpha_{y} K(x,c(Q))\right| & \\
 & \hspace{-2cm} \lesssim \sum_{|\alpha|=N_{q_1}+1} \left\|\frac{(y-c(Q))^\alpha}{\alpha!} \partial^\alpha_{y} K(x,y)\right\|_{\infty,y\in Q}  \\
 &  \hspace{-2cm} \lesssim  |Q|^{(N_{q_1}+1)/d} \frac{h(x)}{|x-c(Q)|^{d+N_{q_1}+1}}. 
\end{align*}
We also get
\begin{align*}
 \left|T(a)(x)\right| & \lesssim \int |Q|^{(N_{q_1}+1)/d} \frac{h(x)}{|x-c(Q)|^{d+N_{q_1}+1}} |a(y)| dy \\
  & \lesssim |Q|^{(N_{q_1}+1)/d} \frac{h(x)}{|x-c(Q)|^{d+N_{q_1}+1}} |Q|^{1-1/q_1}.
 \end{align*}
Therefore with the Hölder inequality on the weak Lebesgue spaces $L^{p,\infty}$ and by integrating $x\in (5Q)^{c}$, we obtain
\begin{align*}
 \left\|T(a)\right\|_{q,\infty,(5Q)^c} & \lesssim |Q|^{(N_{q_1}+1)/d} \|h\|_{r,\infty} \frac{1}{|Q|^{1-1/q_1+(N_{q_1}+1)/d}} |Q|^{1-1/q_1} \\
 & \lesssim 1.
  \end{align*}
We compute the proof by studying the case $x\in 5Q$ with the Hölder inequality and the $L^{p}$-boundedness of $T$~:
\begin{align*}
 \left\|T(a)\right\|_{q,\infty,(5Q)} & \lesssim |Q|^{1/q-1/p} \left\|T(a)\right\|_{p} \\
 & \lesssim |Q|^{1/q-1/p} \|a\|_{p_1} \lesssim |Q|^{1/q-1/p + 1/p_1-1/2+1/2-1/{q_1}} \lesssim 1.
\end{align*}

\mb Here we have assumed that $p_1\leq 2$. If it is not the case, we have to consider the $L^s$-atoms of $H^{q_1}$ with $s\geq p_1$ or use first our previous theorem to have continuities for $T$ with $p_1=1$. By consequence, we have shown that $T$ is bounded on all the atoms of $H^{q_1}$ into $L^{q,\infty}$. 
\findem

\begin{rem} Nowadays, it is well known that an operator, which is bounded on whole the set of atoms, does not always admit a continuous extension to the whole Hardy space.
There is a counterexample for the Hardy space $H^1$ in \cite{meyer}.
\end{rem}

\mb We now apply this abstract result to our operator.

\begin{prop} \label{para3}  Let $f_1,..,f_{n-2}$ be fixed and smooth functions belonging to $\s(\R^d) \subset L^\infty$ and $f_{n}$ be a smooth function belonging to $H^{q}$  with $q\leq 2$ (or $L^{1}$). Then the operator
$$V := f_{n-1} \rightarrow U_{\rho,\lambda,L}(f_1,..,f_{n})$$
satisfies the assumptions of Theorems \ref{L1faible} and \ref{hardy} for $p_1=\infty$, $p=q=r$ and $h \lesssim M^{**}_{d/q+1}(f_{n},\Phi^{n})$ at any order $N$. In addition the bounds can be uniformly controlled with respect to $\rho$ and $\lambda$ under the condition (\ref{conditionmu}). Here we set $M^{**}_{b}$ for the following maximal operator~:
$$M^{**}_{b}(f,\Phi^n)(x) := \sup_{t>0} \sup_{y\in \R^d} \left(1+t^{-1}|y|\right)^{-b} \left|\Phi^n_t \ast f(x-y) \right|.$$
\end{prop}

\dem The assumption of the boundedness is given by Corollary \ref{cor1}. So we have just to check the assumption about the kernel $K(x,z)$, which is given by
\begin{eqnarray*}
\lefteqn{ K(x,z) =\int_0^\infty \int_{\R^d}  \Psi_{t}(y) \Phi^{n-1}_{\lambda_{n-1}t}(x-\rho_{n-1}\lambda_{n-1}y-z) } & & \\
& & \prod_{j=1}^{n-2} \Phi^{j}_{\lambda_{j}t} \ast f_j (x-\rho_{j}\lambda_{j}y)
   \Phi^{n}_{\lambda_{n}t} \ast f_{n} (x-\rho_{n}\lambda_{n}y) \frac{L(t)}{t} dydt.
\end{eqnarray*}
We can differentiate the kernel and we obtain
\begin{align*}
\lefteqn{ \left| \partial^\beta_{z} K(x,z)\right| \lesssim \int_0^\infty
\prod_{j=1}^{n-2} \|f_j\|_\infty  \int_{\R^{d}} \left(1+\frac{|x-\rho_{n-1}\lambda_{n-1}y-z|}{\lambda_{n-1}
t}\right)^{-M} } & & \\
& &   |\lambda_{n-1}t|^{-d-|\beta|}  \left(1+\frac{|y|}{t}\right)^{-2M}\left(1+\frac{|y|}{t}\right)^{-d/q+1} M^{**}_{d/q+1}(f_{n},\Phi^{n})(x)  \frac{\|L\|_\infty dydt}{t^{d+1}}.
\end{align*}
With the arguments, used in the proof of Proposition \ref{propnoyau}, we can estimate the integrals. We also get the desired result~:
$$ \left|\partial_z^\beta K(x,z)\right| \lesssim \prod_{j=1}^{n-2} \|f_j\|_\infty \|L\|_\infty M^{**}_{d/q+1}(f_{n},\Phi^{n})(x) |x-z|^{-d-|\beta|},$$
uniformly with respect to $\lambda$. \findem

\mb
\begin{cor} \label{cor2} The operator $U_{\rho,\lambda,L}$ can be continuously extended from $(L^\infty_c)^{\otimes(n-2)} \times H^{p} \times H^{q}$ into $L^{s}$ for all exponents $0<q\leq 2$ and $0< p \leq \infty$ such that
$$\frac{1}{s}=\frac{1}{p}+\frac{1}{q}.$$
In addition if $p=1$ or $q=1$, we are allowed to substitute the Hardy space $H^{1}$ by the Lebesgue space $L^{1}$ with changing the final space $L^{s,\infty}$ instead of $L^{s}$.
All these continuity bounds are uniform on $0<\rho_i\leq 1$ and $\lambda$ satisfying (\ref{conditionmu}).
\end{cor}

\dem It is a direct consequence of the previous Proposition and the two previous Theorems. We use the fact the maximal operator $M^{**}_{d/q+1}$ is continuous from $H^{q}$ to $L^{q}$ for all exponent $0<q<\infty$ and from $L^1$ to $L^{1,\infty}$. This claim is proved in Theorem 6.4.4 of \cite{Gra}.
So for $f_1,..,f_{n-2}$ fixed bounded and compactly supported functions and $f_{n}\in H^q$, we obtain that the operator 
$$V:=f_{n-1} \rightarrow U_{\rho,\lambda,L}(f_1,..,f_{n})$$
is bounded on all the $H^p$-atoms into $L^s$. Now we use that $U_{\rho,\lambda,L}$ is bounded from $(L^2)^{\otimes (n-2)} \times H^p \times H^q$ into $L^t$ (for the corresponding exponent $t$, see Theorem \ref{thm:final1}) and that the functions $f_i$ are in $L^2$ (beeing compactly supported) in order to be able to extend $V$ on the whole Hardy space $H^p$. This is a classical argument (see for example the proof of Theorem 6.7.1 in \cite{Gra}). We also obtain the continuity of $U_{\rho,\lambda,L}$ from $(L^\infty_c)^{\otimes(n-2)} \times H^{p} \times H^{q}$ into $L^{s}$. We use the same ideas for the space $L^1$ instead of $H^p$ with $p=1$.
\findem

\mb By producing the same reasoning over each component and by using interpolation results, we can prove Theorem \ref{thm:general} for our model operators.

\mb {\noindent {\bf Proof of Theorem \ref{thm:general} for the model operators : }} The case where $S_2=\emptyset$ was shown in Section \ref{section:littlewood} : Theorem \ref{thm:final1} (with the discrete equivalent model) and is a consequence of Theorem \ref{L1faible} and Proposition \ref{para2} for the continuities in weak Lebesgue spaces. The case where $S_2\neq \emptyset$ is a consequence of Theorems \ref{th7}, \ref{L1faible}, \ref{hardy} and Proposition \ref{para2}.  \findem

\section{Decomposition of multipliers with our model operators}
\label{section:paradecomp}

In this section we will prove how to decompose a multilinear multiplier of Theorem \ref{thm:general} with our model operators. This reduction will also conclude the proof of this Theorem. The way to decompose a multilinear multiplier with paraproducts is well known (see for example \cite{cm}). We quickly remember this operation and check that we keep the uniformity with respect to the important parameter $\lambda$.

\gb So let $T$ be an operator of Theorem \ref{thm:general}. It is also associated to a symbol $\sigma$ which satisfies
\be{condising} \forall m_i \in \N^{d} \qquad \left|\partial_{\xi_1}^{m_1} .. \partial_{\xi_{n-1}}^{m_{n}} \sigma (\xi_1,..,\xi_{n}) \right| \lesssim \frac{\prod_{i=1}^{n} |\lambda_i^{|m_i|}|} {d_\lambda(\xi,0)^{|m_1|+..+|m_{n-1}|}}. \ee
As we have seen in Proposition \ref{proppara}, our model operators allow us to get the paraproducts. So we use the ``classical'' decomposition of an H\"ormander multiplier with paraproducts. Let us recall it (we use the ideas of \cite{cm}). \\

\mb For any index $l\in\{1,..,n\}$, we choose a smooth homogeneous function $\zeta_l$ on $(\R^d)^{n}$ supported in the cone~:
$$\left\{ \xi\in (\R^d)^{n},\  |\lambda_l \xi_l|\simeq \max_{j} |\lambda_j \xi_j| \right\}.$$
We can choose them in order that
$$\forall \xi\in (\R^d)^{n}, \qquad \sum_{l=1}^{n} \zeta_{l}(\xi) =1.$$
Let $\Psi$ be a real and smooth function on $\R^d$ whose the spectrum is contained in a corona around $0$ and such that
\be{conditions2} \forall \eta\in \R^d\setminus\{0\} \qquad
\sum_{k\in\Z^d} \left|\widehat\Psi(2^{k} \eta)\right|^2=1.\ee
Let $\phi$ be a smooth function on $\R^d$ whose the spectrum is bounded and such that
$$\forall l\in\{1,..,n\},\ \forall \xi\in supp(\zeta_l) \qquad \widehat{\Psi}(2^k\lambda_l\xi_{l}) \neq 0 \Longrightarrow \ \forall j\neq l,\  \widehat{\Phi}(2^{k} \lambda_j\xi_j) =1.$$
We have also a partition of the unity on the whole frequency plane and we get
\begin{eqnarray*}
\lefteqn{T(f_1,..,f_{n})(x)=} \\
 & & \sum_{k\in\Z} \sum_{l=1}^{n} \int_{\R^{dn}} e^{ix.(\xi_1+..+\xi_{n})}
\widehat{\Psi}^2(2^k\lambda_l \xi_l)
\widehat{f_{l}}(\xi_{l})\prod_{j\neq l} \widehat{\Phi}^2(2^k\lambda_j \xi_j)\widehat{f_j}(\xi_j) \sigma(\xi) d\xi.
\end{eqnarray*}
Let us define new symbols
$$\sigma_{l,k}(\xi):=\sigma(\frac{\xi_1}{\lambda_1 2^k},..,\frac{\xi_{n}}{\lambda_{n} 2^k}) \widehat{\Psi}(\xi_l)\prod_{j\neq l} \widehat{\Phi}(\xi_j).$$
Hence the operator $T$ can be written by
\begin{align*}
\lefteqn{T(f_1,..,f_{n})(x) = \sum_{k\in\Z } \sum_{l=1}^{n} \int_{\R^{dn}} e^{ix.(\xi_1+..+\xi_{n})} } & & \\
& &  \widehat{\Psi}(2^k\lambda_l \xi_l)
\widehat{f_{l}}(\xi_{l})\prod_{j\neq l} \widehat{\Phi}(2^p\lambda_j \xi_j)\widehat{f_j}(\xi_j) \sigma_{l,k}(2^{k}\lambda_1\xi_1,..,2^k\lambda_{n}\xi_{n}) d\xi.
\end{align*}
With the assumption (\ref{condising}), we get that
$$\sigma_{l,k} \in L^1 \qquad \textrm{and} \qquad  \Delta^{N}\sigma_{l,k} \in L^1$$
uniformly with respect to $k$ and $l$ for an integer $N$ as large as we want. So the symbols $\sigma_{l,k}$ satisfy
$$\sigma_{l,k}(\xi)=\int_{\R^{dn}} e^{i\xi.u} \frac{L(l,k,u)}{\left(1+|u|^2\right)^{N}} du$$
with a function $L\in L^\infty(\{1,..,n\} \times \Z \times \R^{dn})$. Then we have 
\begin{align*}
T(f_1,..,f_{n})(x) & =  \sum_{k\in\Z} \sum_{l=1}^{n} \iint_{\R^{2dn}} e^{ix.(\xi_1+..+\xi_{n})} \frac{L(l,k,u)}{\left(1+|u|^2\right)^{N}}  \\
 & \quad \widehat{\Psi}(2^k\lambda_l \xi_l)
\widehat{f_{l}}(\xi_{l})\prod_{j\neq l} \widehat{\Phi}(2^k\lambda_j \xi_j)\widehat{f_j}(\xi_j) e^{i2^{-k}(\lambda_1^{-1}\xi_1,..,\lambda_{n}^{-1}\xi_{n}).u} d\xi du \\
& = \int_{\R^{dn}} \sum_{k} \sum_{l=1}^{n} \frac{L(l,k,u)}{\left(1+|u|^2\right)^{N}} \int_{\R^{dn}} e^{ix.(\xi_1+..+\xi_{n})} \\ 
& \quad  \widehat{\tau_{u_l}\Psi}(2^k\lambda_l \xi_l)
\widehat{f_{l}}(\xi_{l})\prod_{j\neq l} \widehat{\tau_{u_j}\Phi}(2^k\lambda_j \xi_j)\widehat{f_j}(\xi_j) d\xi du.
\end{align*}
Here we are writing $\tau_{y}$ for the translation of the vector $y\in\R^d$. We also obtain
\begin{align} 
\lefteqn{T(f_1,..,f_{n-1})(x) = \int_{\R^{dn}} \sum_{l=1}^{n} \frac{1}{\left(1+|u|^2\right)^{N}} } \nonumber & & \\
 & &   \sum_{k\in\Z}  L(l,k,u)
\left[(\tau_{u_l}\Psi)_{\lambda_l 2^{k}}\ast f_{l}(x)\right] \prod_{j\neq l} \left[(\tau_{u_j}\Phi)_{2^{k}\lambda_j} \ast f_j(x) \right] du. \label{decompos} 
\end{align}
For $l$ and $u$ being fixed, we also find ``classical'' paraproducts.

\mb We now can finish the proof of Theorem \ref{thm:general}~:

\gb {\noindent {\bf End of the proof of Theorem \ref{thm:general} : }} \\
We have seen in Proposition \ref{proppara} that the ``classical'' paraproducts are obtained as limit objects of our models operators when some $\rho_i$ tends to $0$.
So the uniform results of Theorem \ref{thm:general}, proved for our model operators (at the end of the previous section), are also satisfied for the paraproducts. \\
So for each $l$ and $u$ fixed, the operator
$$ (f_1,..,f_n) \rightarrow \sum_{k\in\Z} L(l,k,u)
\left[(\tau_{u_l}\Psi)_{\lambda_l 2^{k}}\ast f_{l}\right] \prod_{j\neq l} \left[(\tau_{u_j}\Phi)_{2^{k}\lambda_j} \ast f_j \right] $$
satisfies all the continuities of Theorem \ref{thm:general}. These continuities are bounded by a weight $(1+|u|)^K$ for a large enough integer $K$ (uniformly with respect to $\lambda$). \\
So if the exponent of the final space $p$ is bigger than $1$, by using the triangular inequality with an integer $N\gg K$, we get the same continuities for the operator $T$. \\ 
If $p<1$  we cannot use the triangular inequality. \\
$*$ If all the exponents are finite (first part with the Littlewood-Paley square functions) : we exactly use the same proof. The spectral study is identical due to the fact that the parameters $u$ and $l$ have no importance. With Lemma \ref{ponctuel} we can have a pointwise estimates of the different square functions which permit us to obtain the result.\\
$*$ If some exponent are infinite (second part with Carleson measure and Calder\'on-Zygmund decomposition), the proof is based on the first continuity from $(L^\infty)^{n-1} \times L^2$ into $L^2$ (which is satisfied for $T$ by the triangular inequality) and on estimates about the multilinear kernel (which are again satisfied for $T$ by the linear correspondance between the kernel and the operator). \\
In the two cases, the continuities of Theorem \ref{thm:general} are proved for $T$.
 \findem

\mb We have also finish the proof of our Theorem \ref{thm:general}. A question stays open : is one of our condition a-) or b-) (in Theorem \ref{thm:general})  necessary to have uniform estimates ?


\begin{thebibliography}{10}

\bibitem{bony}
J.M. Bony,
{{C}alcul symbolique et propagation des singularit\'es pour les
  \'equations aux d\'eriv\'ees partielles non lin\'eaires},
\emph{Ann. Sci. Eco. Norm. Sup.} {\bf 14}, 209-246, (1981).

\bibitem{Bownik}
 M. Bownik,
{ Boundedness of operators on Hardy spaces via atomic decompositions}, 
\emph{Proc. Amer. Math. Soc.} \textbf{133} (12), 3535-3542, (2007).

\bibitem{CM2}
R.~Coifman and Y.~Meyer,
{ On commutators of singular integrals and bilinear singular integrals},
\emph{Trans. Amer. Math. Soc.} \textbf{212}, 315-331, (1975).

\bibitem{cm}
R.~Coifman and Y.~Meyer,
\emph{ Au del\`a des op\'erateurs pseudo-diff\'erentiels},
\newblock SMF Astérisque \textbf{57}, (1978).

\bibitem{CM3}
R.~Coifman and Y.~Meyer.
{Commutateurs d'int\'egrales singuli\`eres et op\'erateurs multilin\'eaires},
\emph{Ann. Inst. Grenoble} \textbf{28}, 177-202, (1978).

\bibitem{CM4}
R.~Coifman and Y.~Meyer,
\emph{Ondelettes et op\'erateurs III, O\'erateurs multilin\'eaires},
\newblock Actualit\'s Math\'ematiques, Hermann, (1991).

\bibitem{lifan}
D.~Fan and X.~Li
{A bilinear oscillatory integrals along parabolas},
\emph{ Preprint, arXiv:0709.2907}.

\bibitem{Gra}
L.~Grafakos,
\emph{ Classical and Modern Fourier Analysis},
\newblock Pearson Education, (2004).

\bibitem{GK}
L.~Grafakos and N.~Kalton,
{The {M}arcinkiewicz multiplier condition for bilinear operators},
\emph{Studia Math.} \textbf{146} no.2, 115-156, (2001).

\bibitem{GK2}
L.~Grafakos and N.~Kalton,
{Multilinear Calder\'on-Zygmund operators on {H}ardy spaces},
\emph{Coll. Math.} \textbf{52}, 169-179, (2001).


\bibitem{GT}
L.~Grafakos and R.~Torres,
{Multilinear {C}alder\'on-{Z}ygmund theory},
\emph{Adv. in Math.} \textbf{165}, 124-164, (2002).

\bibitem{KS}
C.~Kenig and E.M. Stein,
{Multilinear estimates and fractional integration},
\emph{Math. Res. Lett.} \textbf{6} no.1, 1-15, (1999).

\bibitem{li}
X.~Li,
{Uniform estimates for some paraproducts},
\emph{Preprint, arXiv : 0709.2906}.

\bibitem{MSV}
S. Meda, P. Sj\"ogre and M. Vallarino, 
{On the $H^1-L^1$ boundedness of operators},
\emph{Proc. Amer. Math. Soc.}, (2008).

\bibitem{meyer}
Y. Meyer, M. Taibleson and G. Weiss, 
{Some functional analytic properties of the spaces $B_q$ generated by blocks}, 
\emph{Indiana. Univ. Math. J.} \textbf{34}, 493-515, (1985).

\bibitem{mptt}
C.~Muscalu, J.~Pipher, T.~Tao, and C.~Thiele,
{A short proof of the {C}oifman-{M}eyer multilinear theorem},
\newblock Non published, \emph{available at
  http://www.math.brown.edu/~jpipher/trilogy1.pdf} .

\bibitem{MTT4}
C.~Muscalu, T.~Tao, and C.~Thiele,
{Uniform estimates on paraproducts},
\emph{Journ. Anal. Math} \textbf{87}, 369-384, (2002).

\bibitem{mtt}
C.~Muscalu, T.~Tao, and C.~Thiele,
{Uniforms estimates on multi-linear operators with modulation
  symmetry},
\emph{Journ. Anal. Math} \textbf{88}, 255-309, (2002).

\bibitem{S}
E.M. Stein,
\emph{Singular integrals and differentiability properties of
  functions},
\newblock Princeton Univ. Press, (1970).

\bibitem{stein}
E.M. Stein,
\emph{Harmonic analysis : Real variable Methods, Orthogonality, and
  Oscillatory Integrals},
\newblock Princeton Univ. Press, (1993).

\bibitem{Uchi}
A.~Uchiyama,
\emph{Hardy Spaces on the Euclidean Space},
\newblock Springer, (2001).

\end{thebibliography}
\end{document}